\documentclass[a4paper,12pt,reqno]{amsart}
\usepackage{amsmath}
\usepackage{amssymb}
\usepackage{amsthm}

\usepackage{amscd}
\usepackage{amsfonts}
\usepackage{setspace}
\usepackage{version}

\title{A note on character sums over short moving intervals}
\author{Adam J Harper}
\address{Mathematics Institute, Zeeman Building, University of Warwick, Coventry CV4 7AL, England}
\email{A.Harper@warwick.ac.uk}
\date{17th March 2022}
\thanks{Quite a lot of the research leading to this paper was carried out around 2016, when the author was supported by a research fellowship at Jesus College, Cambridge.}

\numberwithin{equation}{section}
\theoremstyle{plain}

\onehalfspace
\newcommand{\N}{\mathbb{N}}
\newcommand{\R}{\mathbb{R}}
\newcommand{\E}{\mathbb{E}}

\newcommand{\C}{\mathbb{C}}

\newtheorem{conj1}{Conjecture}
\newtheorem{thm1}{Theorem}
\newtheorem{thm2}[thm1]{Theorem}
\newtheorem{thm3}[thm1]{Theorem}
\newtheorem{thm4}[thm1]{Theorem}

\newtheorem{prob1}{Probability Result}
\newtheorem{numb1}{Number Theory Result}
\newtheorem{numb2}[numb1]{Number Theory Result}

\newtheorem{prob2}[prob1]{Probability Result}
\newtheorem{prob3}[prob1]{Probability Result}
\newtheorem{cl1}{Claim}
\newtheorem{cl2}[cl1]{Claim}
\newtheorem{cl3}[cl1]{Claim}
\newtheorem{numb3}[numb1]{Number Theory Result}

\begin{document}

\maketitle

\begin{abstract}
We investigate the sums $(1/\sqrt{H}) \sum_{X < n \leq X+H} \chi(n)$, where $\chi$ is a fixed non-principal Dirichlet character modulo a prime $q$, and $0 \leq X \leq q-1$ is uniformly random. Davenport and Erd\H{o}s, and more recently Lamzouri, proved central limit theorems for these sums provided $H \rightarrow \infty$ and $(\log H)/\log q \rightarrow 0$ as $q \rightarrow \infty$, and Lamzouri conjectured these should hold subject to the much weaker upper bound $H=o(q/\log q)$. We prove this is false for some $\chi$, even when $H = q/\log^{A}q$ for any fixed $A > 0$. On the other hand, we show it is true for ``almost all'' characters on the range $q^{1-o(1)} \leq H = o(q)$.

Using P\'{o}lya's Fourier expansion, these results may be reformulated as statements about the distribution of certain Fourier series with number theoretic coefficients. Tools used in the proofs include the existence of characters with large partial sums on short initial segments, and moment estimates for trigonometric polynomials with random multiplicative coefficients.
\end{abstract}

% INTRODUCTION %%%%%%%%%%%%%%%%%%%%%%%%%%%%%
\section{Introduction}
Let $q$ denote a large prime, and $\chi$ a non-principal Dirichlet character modulo $q$. In this paper we will be interested in the statistical behaviour of sums
$$ S_{\chi,H}(x) := \sum_{x < n \leq x+H} \chi(n) , $$
where $H = H(q)$ is some function.

Since $\chi$ has period $q$, we may restrict attention to $1 \leq H \leq q$. The case of long sums, where $H(q) \asymp q$ as $q \rightarrow \infty$, has been quite extensively studied. See, for example, the work of Granville and Soundararajan~\cite{gransoundpv} and of Bober, Goldmakher, Granville and Koukoulopoulos~\cite{bobgoldgrankou} investigating the largest possible values of character sums, and the recent work of Hussain~\cite{hussain} on the behaviour of the paths $t \mapsto \sum_{n \leq qt} \chi(n)$. In this paper we focus instead on short sums, where $H(q) = o(q)$ as $q \rightarrow \infty$. Our primary focus shall be on the situation where $\chi$ is fixed for given $q$, and the start point $x \in \{0,1,...,q-1\}$ varies, although we will touch on what happens when $\chi$ varies as well.

This problem was studied by Davenport and Erd\H{o}s~\cite{davenporterdos}, who proved that if $\chi = \left(\frac{\cdot}{q}\right)$ is the Legendre symbol; and if the function $H$ satisfies $H \rightarrow \infty$ but $(\log H)/\log q \rightarrow 0$ as $q \rightarrow \infty$; and if $X \in \{0,1,...,q-1\}$ is uniformly random; then one has convergence in distribution to a standard Gaussian,
$$ \frac{S_{\chi,H}(X)}{\sqrt{H}} \stackrel{d}{\rightarrow} N(0,1) \;\;\;\;\; \text{as} \; q \rightarrow \infty . $$
Lamzouri~\cite{lamzouridcs} recently extended this to more general Dirichlet characters. He showed that if one chooses a non-real character $\chi$ modulo each prime $q$ (in any way), then under the same conditions on $H$ as Davenport and Erd\H{o}s~\cite{davenporterdos} one has
$$ \frac{S_{\chi,H}(X)}{\sqrt{H}} \stackrel{d}{\rightarrow} Z_1 + iZ_2 \;\;\;\;\; \text{as} \; q \rightarrow \infty , $$
where $Z_1 , Z_2$ are independent $N(0,1/2)$ random variables. Lamzouri~\cite{lamzouridcs} also obtained a quantitative rate of convergence (in the sense of Kolmogorov distance). We also mention slightly earlier work of Mak and Zaharescu~\cite{makzah}, who proved separate distributional convergence results for the real and imaginary parts of $\frac{S_{\chi,H}(X)}{\sqrt{H}}$, and more generally for the projections of various kinds of moving character sum onto lines through the origin.

All of these results, and many related ones (e.g. the work of Perret-Gentil~\cite{perret} on short sums of $l$-adic trace functions), ultimately depend on a moment method. For example, in the case $\chi = \left(\frac{\cdot}{q}\right)$ Davenport and Erd\H{o}s calculated
$$ \frac{1}{q} \sum_{0 \leq x \leq q-1} \left(\frac{S_{\chi,H}(x)}{\sqrt{H}} \right)^j = \frac{1}{q H^{j/2}} \sum_{1 \leq h_1 , ..., h_j \leq H} \sum_{0 \leq x \leq q-1} \left(\frac{x + h_1}{q}\right) \left(\frac{x + h_2}{q}\right) ... \left(\frac{x + h_j}{q}\right) , $$
showing that for each fixed $j \in \N$ this converges to the standard normal moment $(1/\sqrt{2\pi}) \int_{-\infty}^{\infty} z^j e^{-z^{2}/2} dz$ as $q \rightarrow \infty$. It is well known that the normal distribution is sufficiently nice that moment convergence implies convergence in distribution. The key to performing the moment calculation is that for a given tuple $(h_1 , ..., h_j)$ of shifts, if any shift $h$ occurs with odd multiplicity then the sum over $x$ is $\ll_{j} \sqrt{q}$, by the Weil bound. Under the condition $(\log H)/\log q \rightarrow 0$, all these terms together give a contribution
$$ \ll_{j} \frac{1}{\sqrt{q} H^{j/2}} \sum_{\substack{1 \leq h_1 , ..., h_j \leq H , \\ \text{a shift occurs with odd multiplicity}}} 1 \leq \frac{H^{j/2}}{\sqrt{q}} \rightarrow 0 \;\;\;\;\; \text{as} \; q \rightarrow \infty . $$
If one drops the condition $(\log H)/\log q \rightarrow 0$ then this method seems to break down.

Lamzouri~\cite{lamzouridcs} made the following conjecture about what happens for larger $H$.
\begin{conj1}[Lamzouri, 2013]\label{lamzconj}
Suppose that $H \rightarrow \infty$ but $H = o(q/\log q)$ as the prime $q \rightarrow \infty$. Then if $\chi = \left(\frac{\cdot}{q}\right)$ is the Legendre symbol, and if $X \in \{0,1,...,q-1\}$ is uniformly random, we have
$$ \frac{S_{\chi,H}(X)}{\sqrt{H}} \stackrel{d}{\rightarrow} N(0,1) \;\;\;\;\; \text{as} \; q \rightarrow \infty . $$

If we choose a non-real character $\chi$ modulo each prime $q$ (in any way), then on the same range of $H$ we have
$$ \frac{S_{\chi,H}(X)}{\sqrt{H}} \stackrel{d}{\rightarrow} Z_1 + iZ_2 \;\;\;\;\; \text{as} \; q \rightarrow \infty , $$
where $Z_1 , Z_2$ are independent $N(0,1/2)$ random variables.
\end{conj1}

Our goal here is the further investigation of Lamzouri's conjecture. Prior to this, we briefly explain the origins of the conjecture, and in particular of the condition $H = o(q/\log q)$. For each prime $p$, let $f(p)$ be an independent random variable taking values $\pm 1$ with probability 1/2 each (i.e. a Rademacher random variable), and then for each $n \in \N$ define
$$ f(n) := \prod_{p^{\alpha} ||n} f(p)^{\alpha} , $$
where $p^{\alpha} ||n$ means that $p^{\alpha}$ is the highest power of $p$ that divides $n$. We shall refer to such $f$ as an {\em extended Rademacher random multiplicative function}, and think of $f$ as a random model for the Legendre symbol $\left(\frac{n}{q}\right)$ as $q$ varies. Similarly, to model a complex Dirichlet character $\chi(n)$ we let $f(p)$ be uniformly distributed on the complex unit circle (i.e. Steinhaus random variables), and again define $f(n) := \prod_{p^{\alpha} ||n} f(p)^{\alpha}$, a {\em Steinhaus random multiplicative function}. Chatterjee and Soundararajan~\cite{chatsound} showed that, for a very similar\footnote{Chatterjee and Soundararajan~\cite{chatsound} studied Rademacher random multiplicative functions summed over squarefree numbers, rather than the extended functions whose support isn't restricted to squarefree $n$.} kind of real random function $f$, one has
$$ \frac{\sum_{x < n \leq x+y} f(n)}{\sqrt{\E \left(\sum_{x < n \leq x+y} f(n) \right)^2 }} \stackrel{d}{\rightarrow} N(0,1) \;\;\;\;\; \text{as} \; x \rightarrow \infty $$
provided the interval length $y=y(x)$ satisfies $x^{1/5}\log x \ll y = o(x/\log x)$. Lamzouri's conjectured condition $H = o(q/\log q)$ is analogous to Chatterjee and Soundararajan's upper bound on $y$.

There are at least two issues that need to be understood when considering whether the random multiplicative model is a good one for $S_{\chi,H}(X)$. The first is whether a random multiplicative function captures all of the important structure of a Dirichlet character, which in particular has an additional periodicity property. The second is whether one can infer things about $S_{\chi,H}(X)$, where the function $\chi$ is fixed (for given $q$) and the start point $X$ of the interval randomly varies, from things about $\sum_{x < n \leq x+y} f(n)$ where the interval is fixed (for given $x$) and the function $f$ randomly varies. One might think that, if the latter is a good model for character sums, it would rather be for the case of a fixed interval for given $q$ and randomly varying character $\chi$ mod $q$.

\subsection{Statement of results}
Our main results are negative, showing that Conjecture \ref{lamzconj} is {\em not} correct as stated.
\begin{thm1}
Let $A > 0$ be arbitrary but fixed, and set $H(q) = q/\log^{A}q$. Then as $q$ varies over all large primes, with $\chi = \left(\frac{\cdot}{q}\right)$ denoting the unique corresponding quadratic character, we have
$$ \frac{S_{\chi,H}(X)}{\sqrt{H}} \stackrel{d}{\not\rightarrow} N(0,1) \;\;\;\;\; \text{as} \; q \rightarrow \infty . $$
\end{thm1}

\begin{thm2}
Let $A > 0$ be arbitrary but fixed, and set $H(q) = q/\log^{A}q$. Then as $q$ varies over all large primes, there exists a corresponding sequence of non-real characters $\chi$ modulo $q$ for which
$$ \frac{S_{\chi,H}(X)}{\sqrt{H}} \stackrel{d}{\not\rightarrow} Z_1 + iZ_2 \;\;\;\;\; \text{as} \; q \rightarrow \infty , $$
where $Z_1 , Z_2$ are independent $N(0,1/2)$ random variables.
\end{thm2}

It may not be very illuminating just to say that something does not converge to a specified limit object. In fact, in the real case covered by Theorem 1, we will show that there exists an infinite sequence of primes $q$ along which $\frac{S_{\chi,H}(X)}{\sqrt{H}}$ has properties that forbid it from closely approaching the $N(0,1)$ limit. This special sequence consists of primes $q$ for which $\left(\frac{\cdot}{q}\right)$ is ``highly biased'', in the sense that its partial sums up to about $q/H$ are not small. Similarly, in the non-real case covered by Theorem 2, the bad character $\chi$ that we select for each prime $q$ is such that its partial sum up to about $q/H$ has large modulus.

To explain further, if $\chi$ is primitive mod $q$ (so for $q$ prime any non-principal character is admissible), then P\'{o}lya's Fourier expansion for character sums implies that
$$ S_{\chi, H}(x) = \frac{\tau(\chi)}{2\pi i} \sum_{0 < |k| < q/2} \frac{\overline{\chi}(-k)}{k} e(kx/q) (e(kH/q) - 1) + O(\log q) . $$
Here $\tau(\chi)$ denotes the Gauss sum, of absolute value $\sqrt{q}$, and $e(\cdot) = e^{2\pi i \cdot}$ denotes the complex exponential. When $|k| \leq q/H$ we have $(1/k) (e(kH/q) - 1) \approx 2\pi i (H/q)$, and it turns out that (on average over $x$) these are essentially the only terms that make a significant contribution, so $S_{\chi, H}(x)/\sqrt{H} \approx (\tau(\chi)\sqrt{H}/q) \sum_{0 < |k| < q/H} \overline{\chi}(-k) e(kx/q)$.

Now if $H = q/\log^{A}q$, and so $q/H = \log^{A}q$, we can find characters $\chi$ for which $|\sum_{0 < k < q/H} \chi(k)| \gg_{A} q/H$. For such characters, we can think of $S_{\chi, H}(x)/\sqrt{H}$ as having a significant piece resembling the scaled Dirichlet kernel $(\tau(\chi)\sqrt{H}/q) \sum_{0 < |k| < q/H} e(kx/q)$. The Dirichlet kernel certainly does not have Gaussian behaviour as $x$ varies and $q \rightarrow \infty$, in fact (since it has relatively small $L^{1}$ norm) it converges to 0 in probability, which suggests it is unlikely that $S_{\chi, H}(x)/\sqrt{H}$ can converge to the desired Gaussian. This argument can be made rigorous by subtracting a suitable multiple of the Dirichlet kernel from $S_{\chi, H}(x)/\sqrt{H}$, which makes no difference to the putative convergence in distribution but reduces the variance of the sum.

Note that the use of P\'{o}lya's Fourier expansion imports information about the periodicity of $\chi$ mod $q$ into our analysis.

\vspace{12pt}
The characters used in the proofs of Theorems 1 and 2 are quite special, suggesting that Lamzouri's conjecture might be true for {\em almost all} $q$ for real characters, or for almost all characters for each $q$ for non-real characters. Another reason for believing this comes from thinking more carefully about the representation $S_{\chi, H}(x)/\sqrt{H} \approx (\tau(\chi)\sqrt{H}/q) \sum_{0 < |k| < q/H} \overline{\chi}(-k) e(kx/q)$. In a famous classical paper, Salem and Zygmund~\cite{salzyg} showed that for almost all sequences of independent $\pm 1$ coefficients, the partial Fourier series with those coefficients satisfy a central limit theorem when the ``frequency'' (corresponding to $x/q$ in our setup) is chosen uniformly at random. Thus, if we believe that the values of a typical Dirichlet character are somewhat ``random looking'', we might expect to have a central limit theorem as the length $q/H$ tends to infinity. This translates into a condition $H = o(q)$, rather than the condition $H = o(q/\log q)$ proposed by Lamzouri~\cite{lamzouridcs}.

In this positive ``almost all'' direction, we establish the following.

\begin{thm3}
Let $H = H(q)$ satisfy $\frac{\log(q/H)}{\log q} \rightarrow 0$ and $H = o(q)$ as the prime $q \rightarrow \infty$. Then there exists a subset $\mathcal{P}_{H}$ of primes, which satisfies $\frac{\#(\mathcal{P}_{H} \cap [Q,2Q])}{\#\{ Q \leq q \leq 2Q : \; q \; \text{prime}\}} \geq 1 - O(e^{-\min_{Q \leq q \leq 2Q}\log^{3/4}(q/H)})$ (say) for all $Q = 2^j, j \in \N$, such that if $\chi = \left(\frac{\cdot}{q}\right)$ we have
$$ \frac{S_{\chi,H}(X)}{\sqrt{H}} \stackrel{d}{\rightarrow} N(0,1) \;\;\;\;\; \text{as} \; q \rightarrow \infty , \; q \in \mathcal{P}_{H} . $$
\end{thm3}

\begin{thm4}
Let $H = H(q)$ satisfy $\frac{\log(q/H)}{\log q} \rightarrow 0$ and $H = o(q)$ as the prime $q \rightarrow \infty$. Then there exist sets $\mathcal{G}_{q, H}$ of characters mod $q$, satisfying $\#\mathcal{G}_{q,H} \geq q(1 - O(e^{-\log^{3/4}(q/H)}))$, such that for any choice of $\chi \in \mathcal{G}_{q, H}$ we have
$$ \frac{S_{\chi,H}(X)}{\sqrt{H}} \stackrel{d}{\rightarrow} Z_1 + iZ_2 \;\;\;\;\; \text{as} \; q \rightarrow \infty , $$
where $Z_1 , Z_2$ are independent $N(0,1/2)$ random variables.
\end{thm4}

The proofs of Theorems 3 and 4 again use the trigonometric series approximation to $S_{\chi, H}(x)/\sqrt{H}$, which can be reworked slightly into a form (roughly speaking) like $S_{\chi, H}(x)/\sqrt{H} \approx \frac{2\sqrt{q}}{\pi \sqrt{H}} \sum_{0 < k < q/H} \frac{\overline{\chi}(k) \sin(\pi k H/q)}{k} \cos(2 \pi k x/q)$. In fact, looking at $S_{\chi, H}(X)/\sqrt{H}$ for $X \in \{0,1,...,q-1\}$ uniformly random turns out to be roughly equivalent to looking at $\frac{2\sqrt{q}}{\pi \sqrt{H}} \sum_{0 < k < q/H} \frac{\overline{\chi}(k) \sin(\pi k H/q)}{k} \cos(2 \pi k \theta)$, for $\theta \in [0,1]$ uniformly random. This latter small change is not really important, but neatens the writing.

Since moment convergence implies distributional convergence to the Gaussian, to prove Theorem 3 it would suffice (roughly speaking) to show the existence of a subsequence $\mathcal{P}_H$ of primes such that, for each fixed $j \in \N$, we have
$$ \int_{0}^{1} \left( \frac{2\sqrt{q}}{\pi \sqrt{H}} \sum_{1 \leq k < q/H} \frac{\left(\frac{k}{q}\right) \sin(\pi k H/q)}{k} \cos(2\pi k\theta) \right)^{j} d\theta \rightarrow (1/\sqrt{2\pi}) \int_{-\infty}^{\infty} z^j e^{-z^{2}/2} dz $$
as $q \rightarrow \infty , q \in \mathcal{P}_{H}$. To do this, we can try to calculate the average (square) discrepancy between the actual and the Gaussian moments as $q$ varies in each dyadic interval, namely
$$ \frac{\log Q}{Q} \sum_{\substack{Q \leq q \leq 2Q, \\ q \; \text{prime}}} \left| \int_{0}^{1} \left( \frac{2\sqrt{q}}{\pi \sqrt{H}} \sum_{1 \leq k < q/H} \frac{\left(\frac{k}{q}\right) \sin(\pi k H/q)}{k} \cos(2\pi k\theta) \right)^{j} d\theta - \frac{1}{\sqrt{2\pi}} \int_{-\infty}^{\infty} z^j e^{-z^{2}/2} dz \right|^2 . $$
If this tends to zero at a sufficient rate as $Q \rightarrow \infty$, on a range of $j$ that grows to infinity as $Q \rightarrow \infty$ (recall that we need convergence of {\em all} fixed integer moments to guarantee convergence to the Gaussian), then we can form $\mathcal{P}_H$ from all the many primes in each interval $[Q,2Q]$ where the discrepancy is simultaneously small for a suitable range of $j$.

Provided that $(q/H)^j$ is small compared with $Q$, so the periodicity of the characters $\left(\frac{k}{q}\right)$ doesn't intervene, one expects the left hand side in the above display to be close to the corresponding one where $\left(\frac{k}{q}\right)$ is replaced by an extended Rademacher random multiplicative function $f(k)$, and the normalised sum $\frac{\log Q}{Q} \sum_{\substack{Q \leq q \leq 2Q, \\ q \; \text{prime}}}$ is replaced by an expectation $\E$. There are technical challenges in establishing this, because $q/H(q)$ might also vary with $q$ in the sum, and averaging over primes $q$ entails non-trivial issues with the distribution of primes, but these problems can be overcome (see Number Theory Result \ref{leglikermf} and section \ref{fixqH} below, essentially one needs to show that the $\left(\frac{k}{q}\right)$ for varying $q$ have similar correlation/orthogonality properties to the random $f(k)$). Unfortunately, the condition that $(q/H)^j$ is small compared with $Q$, for each fixed $j$, forces the unwanted condition $\frac{\log(q/H)}{\log q} \rightarrow 0$ in Theorem 3. This is similar to the condition $(\log H)/\log q \rightarrow 0$ that appeared in the work of Davenport and Erd\H{o}s~\cite{davenporterdos}, Lamzouri~\cite{lamzouridcs}, and others.

Finally we need upper bounds for quantities like
$$ \E \left| \int_{0}^{1} \left( \frac{2\sqrt{Q}}{\pi \sqrt{H}} \sum_{1 \leq k < Q/H} \frac{f(k) \sin(\pi k H/Q)}{k} \cos(2\pi k\theta) \right)^{j} d\theta - \frac{1}{\sqrt{2\pi}} \int_{-\infty}^{\infty} z^j e^{-z^{2}/2} dz \right|^2 , $$
where $f(k)$ is a random multiplicative function. This arithmetic input can be extracted from a nice recent paper of Benatar, Nishry and Rodgers~\cite{bennishrodg}. They were interested in almost sure central limit theorems and size bounds for random trigonometric polynomials $\frac{1}{\sqrt{N}} \sum_{n \leq N} f(n) e(n\theta)$, and directly calculated such expectations using a point counting argument drawing on work of Vaughan and Wooley~\cite{vaughanwooley}. Ultimately one needs to count tuples $(n_1, ..., n_{2j})$ satisfying a small collection of linear and multiplicative equations.

As Benatar, Nishry and Rodgers~\cite{bennishrodg} comment, one can also analyse the distribution of $\frac{1}{\sqrt{N}} \sum_{n \leq N} f(n) e(n\theta)$ using martingale methods, and this was done in unpublished work of the present author (see the paper~\cite{harperlimits} for an application of martingales to a different distributional problem for random multiplicative functions). But to transfer these conclusions to character sums, one would seem to again need moment estimates on the random multiplicative side, not just distributional convergence. These could be obtained (e.g. one can use Burkholder's inequalities~\cite{burkholderaoms} and some calculation to show that all moments remain bounded as $N \rightarrow \infty$, and this combined with distributional convergence implies they must all converge to the desired Gaussian moments), but it seems simpler to rely on the existing calculations of Benatar, Nishry and Rodgers~\cite{bennishrodg}.

In the complex case in Theorem 4, one proceeds exactly similarly in studying the square discrepancy from the moments of the complex Gaussian $Z_1 + iZ_2$, now averaging over all $\chi$ mod $q$ rather than over $Q \leq q \leq 2Q$. Provided that $1 \leq n_1, n_2 < q$, say, we have the identity $\frac{1}{q-1} \sum_{\chi \; \text{mod} \; q} \chi(n_1) \overline{\chi}(n_2) = \textbf{1}_{n_1 \equiv n_2 \; \text{mod} \; q} = \textbf{1}_{n_1 = n_2} = \E f(n_1) \overline{f}(n_2)$, where $f(\cdot)$ denotes a Steinhaus random multiplicative function. This exact equality makes it much easier to establish the connection with random multiplicative functions than in Theorem 3, but the condition $1 \leq n_1, n_2 < q$ ultimately forces the same unwanted constraint $\frac{\log(q/H)}{\log q} \rightarrow 0$.

\subsection{Discussion and open questions}
Our results leave open several problems about the behaviour of $S_{\chi, H}(x)/\sqrt{H}$, and related issues.

The results of Davenport and Erd\H{o}s~\cite{davenporterdos} and of Lamzouri~\cite{lamzouridcs} establish a central limit theorem for all characters provided $H \rightarrow \infty$ but $H = q^{o(1)}$, and our results establish a central limit theorem for almost all characters provided $q^{1-o(1)} \leq H = o(q)$. Moreover, we have shown that one cannot hope to prove a central limit theorem for {\em all} characters when $\log(q/H)/\log\log q$ is bounded. Given this state of affairs, one can ask:
\begin{enumerate}
\item how does $S_{\chi, H}(x)/\sqrt{H}$ behave on the missing range $q^{o(1)} \leq H \leq q^{1-o(1)}$ ?

\item indeed, should it be possible to prove a central limit theorem for {\em all} characters provided $H \rightarrow \infty$ and $\log(q/H)/\log\log q \rightarrow \infty$ ?
\end{enumerate}

The author tentatively conjectures that the answer to (ii) is Yes. In view of Corollary A of Granville and Soundararajan~\cite{gransoundlcs}, if the Generalised Riemann Hypothesis is true then $\sum_{n \leq x} \chi(n) = o(x)$ whenever $\chi$ is a non-principal character mod $q$, and $(\log x)/\log\log q \rightarrow \infty$. This means that, assuming GRH, there would be no construction along the lines of Theorems 1 and 2 available once $\log(q/H)/\log\log q \rightarrow \infty$. So if one believes this is the only barrier to a central limit theorem holding, as is somewhat suggested by the representation $S_{\chi, H}(x)/\sqrt{H} \approx (\tau(\chi)\sqrt{H}/q) \sum_{0 < |k| < q/H} \overline{\chi}(-k) e(kx/q)$ together with the classical work of Salem and Zygmund~\cite{salzyg} on random Fourier series, then one arrives at this conjecture.

However, proving such a result seems difficult. Firstly, the best unconditional estimates we have of the form $\sum_{n \leq x} \chi(n) = o(x)$, where $\chi$ is any non-principal character modulo a prime $q$ (one can sometimes do better for special non-prime moduli), are Burgess-type estimates requiring that $x \geq q^{1/4 - o(1)}$. Thus we would need to assume results like GRH merely to exclude the kind of construction from Theorems 1 and 2 from cropping up. But even allowing such unproved arithmetical results, there is no clear way to go on and establish a central limit theorem on the full range of $H$ in (ii). The problem of understanding the distribution of $(\tau(\chi)\sqrt{H}/q) \sum_{0 < |k| < q/H} \overline{\chi}(-k) e(k\theta)$, where $\theta = X/q$ is random but the coefficients $\overline{\chi}(-k)$ are deterministic, is just one example of the important general problem of understanding the distribution of $\sum_{k} a_k e(k\theta)$, where $a_k$ are interesting deterministic coefficients. See, for example, the work of Hughes and Rudnick~\cite{hughesrudnick} on lattice points in annuli. They encounter similar sums where the $a_k$ involve the number of representations of $k$ as a sum of two squares, and the range of their main theorem involves a similar (conjecturally unnecessary) restriction as in Theorems 3 and 4 to allow a proof by the method of moments.

Indeed, even extending our ``almost all'' results to a wider range of $H$ would be very interesting, and doesn't seem easily attackable.

Another, perhaps rather specialised, question is:
\begin{enumerate}
\setcounter{enumi}{2}

\item what can be said about the distribution of $S_{\chi, H}(X)/\sqrt{H}$, for those characters $\chi$ and interval lengths $H$ where it does not satisfy the expected central limit theorem?
\end{enumerate}

%%Thus the proofs of Theorems 1 and 2 show that (for certain special characters, when $H(q) = q/\log^{A}q$, and after subtracting a suitable piece that tends to zero in probability) the mean square of $S_{\chi, H}(X)/\sqrt{H}$ is too small for it to converge to the target Gaussian distribution. But this doesn't rule out that it could converge in distribution to e.g. a Gaussian with smaller variance, in other words that for such special characters it may be more natural to normalise by a quantity different than $\sqrt{H}$.

\vspace{12pt}
We can also return to the random multiplicative functions $f(n)$ that motivated Lamzouri's conjecture~\cite{lamzouridcs}, and played a role in the proofs of Theorems 3 and 4. As discussed earlier, and perhaps demonstrated by Theorems 1 and 2, the author doesn't believe that Chatterjee and Soundararajan's work~\cite{chatsound} on $\sum_{x < n \leq x+y} f(n)$ provides a natural model for $S_{\chi, H}(X)$. But the study of $\sum_{x < n \leq x+y} f(n)$ is very interesting in its own right. Although Chatterjee and Soundararajan only obtained\footnote{We remark again that Chatterjee and Soundararajan~\cite{chatsound} studied Rademacher random multiplicative functions supported on squarefree numbers only. For the next paragraph, $f(n)$ should be taken to mean this model. Most things discussed will carry over to Steinhaus random multiplicative functions as well, but extended Rademacher random multiplicative functions may exhibit some different behaviour due to significant contributions from squares (on which an extended Rademacher random multiplicative function is identically 1) and numbers with large square factors.} a Gaussian limiting distribution when $x^{1/5}\log x \ll y = o(x/\log x)$, they did not show that their upper bound on $y$ is optimal, and forthcoming work of Soundararajan and Xu~\cite{soundxu} extends the range to $x^{1/5}\log x \ll y \ll \frac{x}{\log^{2\log 2 - 1 + \epsilon}x}$. On the other hand, it follows directly from work of the author~\cite{harperrmflow} that if $\frac{y \sqrt{\log\log x}}{x} \rightarrow \infty$ as $x \rightarrow \infty$, then
$$ \E \frac{|\sum_{x < n \leq x+y} f(n)|}{\sqrt{y}} \leq \frac{\E|\sum_{n \leq x} f(n)| + \E|\sum_{n \leq x+y} f(n)|}{\sqrt{y}} \rightarrow 0 . $$
This implies that $\sum_{x < n \leq x+y} f(n)$ converges in probability to 0, rather than converging to a standard Gaussian, when renormalised by its standard deviation. As Soundararajan and Xu~\cite{soundxu} remark, by looking inside the proofs from \cite{harperrmflow} one can show that $\E \frac{|\sum_{x < n \leq x+y} f(n)|}{\sqrt{y}} \rightarrow 0$ even for somewhat smaller $y$. Thus there is at least one qualitative transition in the distributional behaviour of $\sum_{x < n \leq x+y} f(n)$ when $y$ approaches $x$, and the exact location and nature of this remains to be understood.

As also noted earlier, the author believes that $\sum_{x < n \leq x+y} f(n)$ {\em will} be a good model for the behaviour of $\sum_{x < n \leq x+y} \chi(n)$ where $x, y(x)$ are fixed and the {\em character} $\chi$ varies mod $q$, at least provided $x \leq \sqrt{q}$, say (for $x$ close to $q$, one will again need to be more careful to account for the periodicity of $\chi$). It would be very interesting to obtain rigorous results on the distribution of $\sum_{x < n \leq x+y} \chi(n)$ for varying $\chi$.

Finally, we might wonder:
\begin{enumerate}
\setcounter{enumi}{3}

\item when $f(n)$ is a realisation of a Steinhaus or (extended) Rademacher random multiplicative function, what is the distribution of $\sum_{x < n \leq x+H} f(n)$ as $x$ varies over a long interval?
\end{enumerate}

Although the proofs of Theorems 3 and 4 use the random multiplicative model $f(n)$ for $\chi(n)$, they do not address (iv) because they only use this after first passing to the representation $S_{\chi, H}(x)/\sqrt{H} \approx \frac{2\sqrt{q}}{\pi \sqrt{H}} \sum_{0 < k < q/H} \frac{\overline{\chi}(k) \sin(\pi k H/q)}{k} \cos(2 \pi k x/q)$, the truth of which depends on special properties of Dirichlet characters. Our arguments say nothing directly about the ``model'' object $\sum_{x < n \leq x+H} f(n)$. Of course a little care is required to sensibly interpret question (iv). For example, the function $f(n)$ that is 1 for all $n$ on some long initial segment is a realisation of a random multiplicative function, and has rather exceptional behaviour, but it is a realisation that occurs with extremely small probability. A natural problem might be to investigate the distribution of $\sum_{x < n \leq x+H} f(n)$ for ``most'' realisations of $f$, somewhat analogously to Theorems 3 and 4. Relevant work in the literature includes Najnudel's paper~\cite{najnudelconsec}, which explores the joint distribution of the tuple $(f(x),f(x+1), ..., f(x+H))$ for $x$ varying and $H$ fixed (or slowly growing).

% SECTION 2 %%%%%%%%%%%%%%%%%%%%%%%%%%%%%%
\section{Tools for Theorems 1 and 2}
The proofs of Theorems 1 and 2 rest on the following simple principle.
\begin{prob1}\label{probres1}
Let $0 \leq \tau < 1$, and suppose $(V_n)_{n=1}^{\infty}$ is a sequence of real or complex valued random variables satisfying $\E|V_n|^2 \leq \tau$ for all $n$. Then if $Z$ is any real or complex valued random variable such that $\E|Z|^2 = 1$, we have
$$ V_n \stackrel{d}{\not\rightarrow} Z \;\;\;\;\; \text{as} \; n \rightarrow \infty . $$
\end{prob1}

\begin{proof}[Proof of Probability Result 1]
Choose $a \in \R$ such that $\E\min\{|Z|^2, a^2\} \geq (1+\tau)/2$ (such $a$ exists by the monotone convergence theorem). Since $v \mapsto \min\{|v|^2 , a^2\}$ is a continuous {\em bounded} function on $\C$, if we had $V_n \stackrel{d}{\rightarrow} Z$ then we would have
$$ \E \min\{|V_n|^2 , a^2\} \rightarrow \E \min\{|Z|^2 , a^2\} \;\;\;\;\; \text{as} \; n \rightarrow \infty . $$
But this is impossible, since clearly $\E \min\{|V_n|^2 , a^2\} \leq \E|V_n|^2 \leq \tau < (1+\tau)/2$.
\end{proof}

We remark that although Probability Result 1 is simple, there is a non-trivial issue involved which it is important to recognise. Thus the analogous statement in which for some $\nu > 1$ our sequence satisfied $\E |V_n|^2 \geq \nu$ for all $n$ would be false, as can easily be shown by examples. In general, the failure of moments to converge to the moments of a supposed limit distribution need {\em not}, by itself, imply that convergence in distribution is not happening, since moments may be {\em inflated} by events whose probabilities tend to zero, and which are therefore irrelevant to convergence in distribution.

\vspace{12pt}
As explained in the Introduction, Theorems 1 and 2 will also rely on the existence of non-principal characters with large partial sums. In the non-real case, the existence of such characters follows immediately from work of Granville and Soundararajan~\cite{gransoundlcs}.
\begin{numb1}[see Theorem 3 of Granville and Soundararajan~\cite{gransoundlcs}, 2001]\label{ntgransound}
Let $A > 0$, and suppose $q$ is a prime (say) that is sufficiently large in terms of $A$. Then for any $2 \leq x \leq \log^{A}q$, there exist at least $q^{1 - \frac{2}{\log x}}$ characters $\chi$ mod $q$ for which
$$ \left| \sum_{n \leq x} \chi(n) \right| \geq x \rho(A) \left(1 + O(\frac{1}{\log x} + \frac{\log x (\log\log\log q)^2}{(\log\log q)^2}) \right) , $$
where $\rho(A) > 0$ denotes the Dickman function.
\end{numb1}

In the real case, Granville and Soundararajan (see Theorem 9 of \cite{gransoundlcs}) also proved that for any fixed $A$, if $q$ is large and $x = ((1/3)\log q)^A$ then there exists a fundamental discriminant $q \leq |D| \leq 2q$ for which
$$ \sum_{n \leq x} \left(\frac{D}{n}\right) \geq x (\rho(A) + o(1)) . $$
Unfortunately this isn't quite sufficient for our purposes, because we need to find biased real characters $\left(\frac{\cdot}{q}\right)$ where $q$ is {\em prime}. But by reorganising Granville and Soundararajan's proof a little, and inserting information about the zero-free region and exceptional zeros of Dirichlet $L$-functions, one can prove such a statement. This has been done by Kalmynin~\cite{kalmynin}.
\begin{numb2}[see Theorem 3 of Kalmynin~\cite{kalmynin}, 2019]\label{ntrealchars}
For any fixed $B > 0$, there exists a small constant $c(B) > 0$ such that the following is true. If $Q$ is sufficiently large in terms of $B$, then for any $1 \leq x \leq \log^{B}Q$ there exists a prime $Q < q \leq 2Q$ such that
$$ \sum_{n \leq x} \left(\frac{n}{q}\right) \geq c(B) x . $$
\end{numb2}

\begin{proof}[Proof of Number Theory Result 2]
Theorem 3 of Kalmynin~\cite{kalmynin} directly implies Number Theory Result \ref{ntrealchars} provided that $B \geq 1$ and $x = \log^{B}Q$. However, the lower bound for $S_{0}(Q)$ obtained in Kalmynin's proof implies that one can find primes $Q < q \leq 2Q$ such that $\left(\frac{n}{q}\right) = 1$ for all $n \leq \log^{1/3}Q$. This means that if $x \leq \log^{1/3}Q$, then one can make $\sum_{n \leq x} \left(\frac{n}{q}\right)$ maximally large. And if $\log^{1/3}Q < x \leq \log^{B}Q$, then one can run Kalmynin's proof with the sum over $n \leq \log^{B}Q$ replaced by a sum over $n \leq x$ without changing anything, giving the desired conclusion.
\end{proof}

% SECTION 3 %%%%%%%%%%%%%%%%%%%%%%%%%%%%%%
\section{Proof of Theorem 1}\label{negproofsec}
Let $K \geq 1$. If $\chi$ is a primitive character modulo $q$, then P\'{o}lya's Fourier expansion (see e.g. display (9.19) of Montgomery and Vaughan~\cite{mv}, noting that the restriction $K \leq q^{1-\epsilon}$ there is unnecessary if one is happy with a general error term $\frac{q\log q}{K}$ rather than $\frac{\phi(q) \log q}{K}$) yields that
$$ \sum_{n \leq x} \chi(n) = \frac{\tau(\chi)}{2\pi i} \sum_{0 < |k| \leq K} \frac{\overline{\chi}(-k)}{k} (e(kx/q) - 1) + O(1 + \frac{q\log q}{K}) , $$
so in particular
\begin{equation}\label{polyaq}
S_{\chi, H}(x) = \frac{\tau(\chi)}{2\pi i} \sum_{0 < |k| < q/2} \frac{\overline{\chi}(-k)}{k} e(kx/q) (e(kH/q) - 1) + O(\log q) .
\end{equation}
Here $\tau(\chi)$ denotes the Gauss sum, having absolute value $\sqrt{q}$ for primitive $\chi$.

Recall that $X$ denotes a random variable having the discrete uniform distribution on $\{0,1,...,q-1\}$ (this is the randomness with respect to which we will shortly calculate expectations $\E$), and that $H = H(q) = q/\log^{A}q$, and that $\chi = \left(\frac{\cdot}{q}\right)$ is real-valued in Theorem 1. Next let $0 < \delta \leq 1$ be a parameter, that will be fixed later, and define $\alpha = \alpha(\delta) \in \R$ by
$$ \sum_{1 \leq k \leq \delta q/H} \overline{\chi}(-k) = \sum_{1 \leq k \leq \delta q/H} \chi(-k) = \alpha \sum_{1 \leq k \leq \delta q/H} 1 , $$
and define
$$ G_{\chi, H}(x) := \frac{\alpha \tau(\chi) H}{q} \sum_{1 \leq k \leq \delta q/H} e(kx/q) . $$
As discussed in the Introduction, $G_{\chi, H}(x)$ is the scaled Dirichlet kernel that we shall strategically subtract from $S_{\chi, H}(x)$. (The small parameter $\delta$ is only present for technical reasons, to control lower order terms in Taylor expansions of the complex exponential.)

Before embarking on our main computations, we record some basic observations. By expanding the square and using the fact that $\E e(k_1 X/q) e(-k_2 X/q) = \E e((k_1 - k_2)X/q) = 0$ when $-q/2 < k_1 \neq k_2 < q/2$, we find
$$ \E|G_{\chi, H}(X)|^2 = \frac{\alpha^2 H^2}{q} \E|\sum_{1 \leq k \leq \delta q/H} e(kX/q)|^2 = \frac{\alpha^2 H^2}{q} \sum_{1 \leq k \leq \delta q/H} 1 \leq H , $$
as well as
$$ \E\Biggl|\frac{\tau(\chi)}{2\pi i} \sum_{0 < |k| < q/2} \frac{\overline{\chi}(-k)}{k} e(kX/q) (e(kH/q) - 1) \Biggr|^2 = \frac{q}{(2\pi)^2} \sum_{0 < |k| < q/2} \frac{1}{k^2} |e(kH/q) - 1|^2 . $$
Since $\E|S_{\chi, H}(X)|^2 = \E|\sum_{X < n \leq X+H} \chi(n)|^2 = \sum_{1 \leq h_1, h_2 \leq H} \E \chi(X+h_1) \overline{\chi}(X+h_2) = \sum_{1 \leq h_1, h_2 \leq H} \E \chi(X+h_1 - h_2) \overline{\chi}(X) = H + O(H^{2}/q)$, (explaining why $\frac{S_{\chi, H}(X)}{\sqrt{H}}$ is the natural renormalisation to consider), we deduce that
\begin{eqnarray}
\frac{q}{(2\pi)^2} \sum_{0 < |k| < q/2} \frac{1}{k^2} |e(kH/q) - 1|^2 & = & \E|S_{\chi, H}(X) + O(\log q)|^2 \nonumber \\
& = & \E|S_{\chi, H}(X)|^2 + O(\log q \E|S_{\chi, H}(X)| + \log^{2}q) \nonumber \\
& = & (1+o(1))H \nonumber
\end{eqnarray}
when $H = H(q) = q/\log^{A}q$ (and indeed on a much larger range of $H$ as well).

\vspace{12pt}
Next, using \eqref{polyaq}, expanding the square, and calculating as above (and using the Cauchy--Schwarz inequality and the above estimates of $\E|G_{\chi, H}(X)|^2$ and $\E|S_{\chi, H}(X)|^2$ to control the contribution from the $O(\log q)$ term), we find
\begin{eqnarray}\label{sgdiffdisplay}
&& \E|S_{\chi, H}(X) - G_{\chi, H}(X)|^2 \nonumber \\
& = & \E\Biggl|\tau(\chi) \sum_{1 \leq k \leq \delta q/H} e(kX/q) \left(\frac{\overline{\chi}(-k)}{k} \frac{(e(kH/q) - 1)}{2\pi i} - \frac{\alpha H}{q} \right) + \nonumber \\
&& + \frac{\tau(\chi)}{2\pi i} \sum_{\substack{0 < |k| < q/2 , \\ k \notin [1,\delta q/H]}} \frac{\overline{\chi}(-k)}{k} e(kX/q) (e(kH/q) - 1) + O(\log q) \Biggr|^2 \nonumber \\
& = & q \sum_{1 \leq k \leq \delta q/H} \left|\frac{\overline{\chi}(-k)}{k} \frac{(e(kH/q) - 1)}{2\pi i} - \frac{\alpha H}{q} \right|^2 + \frac{q}{(2\pi)^2} \sum_{\substack{0 < |k| < q/2 , \\ k \notin [1,\delta q/H]}} \frac{1}{k^2} |e(kH/q) - 1|^2 + \nonumber \\
&& + O(\sqrt{H} \log q + \log^{2}q ) .
\end{eqnarray}

Using the Taylor expansion $e(kH/q) = 1 + 2\pi i kH/q + O((kH/q)^2)$, the first sum here is seen to be
\begin{eqnarray}
q \sum_{1 \leq k \leq \delta q/H} \left|\frac{\overline{\chi}(-k)H}{q} + O\left(\frac{\delta H}{q}\right) - \frac{\alpha H}{q} \right|^2 & = & \frac{H^2}{q} \sum_{1 \leq k \leq \delta q/H} \left|\overline{\chi}(-k) + O(\delta) - \alpha \right|^2  \nonumber \\
& = & \frac{H^2}{q} \sum_{1 \leq k \leq \delta q/H} \left|\overline{\chi}(-k) - \alpha \right|^2 + O(\delta^2 H) . \nonumber
\end{eqnarray}
Moreover, since we chose $\alpha$ to be the mean value of $\overline{\chi}(-k) (= \chi(-k))$ over the interval $1 \leq k \leq \delta q/H$, this simplifies to
$$ (1- \alpha^2) \frac{H^2}{q} \sum_{1 \leq k \leq \delta q/H} 1 + O(\delta^2 H) , $$
which we can rewrite (again using the Taylor expansion of the exponential) as
$$ \frac{H^2}{q} \sum_{1 \leq k \leq \delta q/H} 1 - \delta \alpha^2 H + O(\delta^2 H + \frac{H^2}{q}) = \frac{q}{(2\pi)^2} \sum_{1 \leq k \leq \delta q/H} \frac{1}{k^2} |e(kH/q) - 1|^2 - \delta \alpha^2 H + O(\delta^2 H + \frac{H^2}{q}) . $$
Inserting this in \eqref{sgdiffdisplay}, and using our earlier calculation that $\frac{q}{(2\pi)^2} \sum_{0 < |k| < q/2} \frac{1}{k^2} |e(kH/q) - 1|^2 = (1+o(1))H$, we deduce
$$ \E|S_{\chi, H}(X) - G_{\chi, H}(X)|^2 = (1 - \delta \alpha^2 + O(\delta^2) + o(1))H . $$

In particular, note that if $H(q) = q/\log^{A}q$ then we have $\delta q/H \leq q/H = \log^{A}q$. Thus by Number Theory Result \ref{ntrealchars}, there exist arbitrarily large primes $q$ for which, with $\chi = \left(\frac{\cdot}{q}\right)$, we have
$$ \Biggl| \sum_{1 \leq k \leq \delta q/H} \overline{\chi}(-k) \Biggr| = \Biggl| \sum_{1 \leq k \leq \delta q/H} \chi(k) \Biggr| \geq c(A) \delta q/H . $$
In other words, for such $q$ we will have $|\alpha| \geq c(A)$. So if we fix the choice $\delta = c_0 c(A)^2$, for a suitably small absolute constant $c_0 > 0$ to neutralise the implicit constant in the $O(\delta^2)$ term, we will have
\begin{equation}\label{l2bound}
\E\left| \frac{S_{\chi, H}(X) - G_{\chi, H}(X)}{\sqrt{H}} \right|^2 = 1 - \delta\alpha^2 + O(\delta^2) + o(1) \leq 1 - (c_{0}/2) c(A)^4
\end{equation}
whenever $q$ is a large enough prime coming from Number Theory Result \ref{ntrealchars}.

\vspace{12pt}
Now on the other hand, using the formula for summing a geometric progression we may calculate explicitly that, with $|| \cdot ||$ denoting distance to the nearest integer,
$$ \E|G_{\chi, H}(X)| = \frac{|\alpha| H}{\sqrt{q}} \E\left|\sum_{1 \leq k \leq \delta q/H} e(kX/q) \right| \ll \frac{H}{\sqrt{q}} \E \min\{\frac{\delta q}{H}, \frac{1}{||X/q||} \} \ll \frac{H}{\sqrt{q}} (1 + \log(\delta q/H)) , $$
and therefore
$$ \E\left|\frac{G_{\chi, H}(X)}{\sqrt{H}} \right| \ll \frac{(1 + \log(\delta q/H))}{\sqrt{q/H}} \ll \frac{\log(q/H)}{\sqrt{q/H}} \rightarrow 0 \;\;\; \text{as} \; q \rightarrow \infty . $$
By Markov's inequality, it follows that $\frac{G_{\chi, H}(X)}{\sqrt{H}}$ converges in probability to zero as $q \rightarrow \infty$, and so if $\frac{S_{\chi, H}(X)}{\sqrt{H}} \stackrel{d}{\rightarrow} N(0,1)$ then we must also have $\frac{S_{\chi, H}(X) - G_{\chi, H}(X)}{\sqrt{H}} \stackrel{d}{\rightarrow} N(0,1)$.

But combining Probability Result \ref{probres1} with \eqref{l2bound}, we see this convergence in distribution is impossible, which proves Theorem 1.
\qed

% SECTION 4 %%%%%%%%%%%%%%%%%%%%%%%%%%%%%%
\section{Proof of Theorem 2}
The proof of Theorem 2 is extremely similar to that of Theorem 1, so we simply make a few remarks to reassure the reader that no additional difficulties arise.

Indeed, this time we define $\alpha = \alpha(\delta) \in \C$ by
$$ \sum_{1 \leq k \leq \delta q/H} \overline{\chi}(-k) = \alpha \sum_{1 \leq k \leq \delta q/H} 1 , $$
and again we set $G_{\chi, H}(x) := \frac{\alpha \tau(\chi) H}{q} \sum_{1 \leq k \leq \delta q/H} e(kx/q)$. Then the same calculations as in the proof of Theorem 1 show that
\begin{eqnarray}
&& \E|S_{\chi, H}(X) - G_{\chi, H}(X)|^2 \nonumber \\
& = & (1- |\alpha|^2) \frac{H^2}{q} \sum_{1 \leq k \leq \delta q/H} 1 + O(\delta^2 H) + \frac{q}{(2\pi)^2} \sum_{\substack{0 < |k| < q/2 , \\ k \notin [1,\delta q/H]}} \frac{1}{k^2} |e(kH/q) - 1|^2 + \nonumber \\
&& + O(\sqrt{H} \log q + \log^{2}q ) \nonumber \\
& = & (1 - \delta|\alpha|^2 + O(\delta^2) + o(1))H . \nonumber
\end{eqnarray}

Next, in place of Number Theory Result \ref{ntrealchars} we can invoke Number Theory Result \ref{ntgransound}, which implies that for any large prime $q$ we may find a non-real character $\chi$ mod $q$ (in fact several of them) for which
$$ \left| \sum_{1 \leq k \leq \delta q/H} \overline{\chi}(-k) \right| = \left| \sum_{1 \leq k \leq \delta q/H} \chi(k) \right| \geq (\rho(A) + o_{\delta, A}(1)) \sum_{1 \leq k \leq \delta q/H} 1 . $$
For such a character we will have $|\alpha| \geq \rho(A) + o_{\delta, A}(1)$, so if we fix the choice $\delta = c \rho(A)^2$, where $c > 0$ is a suitably small absolute constant, then in place of \eqref{l2bound} we will get
$$ \E\left| \frac{S_{\chi, H}(X) - G_{\chi, H}(X)}{\sqrt{H}} \right|^2 = 1 - \delta|\alpha|^2 + O(\delta^2) + o(1) \leq 1 - (c/2) \rho(A)^4 , $$
provided $q$ is large enough.

Combining this bound with Probability Result \ref{probres1}, and the facts that $\frac{G_{\chi, H}(X)}{\sqrt{H}} \stackrel{p}{\rightarrow} 0$ and that $\E|Z_1 + iZ_2|^2 = \E Z_1^2 + \E Z_2^2 = 1$ (where $Z_1 , Z_2$ are independent $N(0,1/2)$ random variables), we conclude that indeed $\frac{S_{\chi,H}(X)}{\sqrt{H}} \stackrel{d}{\not\rightarrow} Z_1 + iZ_2$ as $q \rightarrow \infty$.
\qed

% SECTION 5 %%%%%%%%%%%%%%%%%%%%%%%%%%%%%%
\section{Tools for Theorems 3 and 4}
As discussed in the Introduction, much of the work in the proofs of Theorems 3 and 4 will be done by some results on random multiplicative functions.

\begin{prob2}[See Theorem 1.1 of Benatar, Nishry and Rodgers~\cite{bennishrodg}]\label{rmfimport}
Let $f(n)$ be an extended Rademacher random multiplicative function. Then uniformly for any large $N$, any coefficients $(a_n)_{n \leq N}$ bounded in absolute value by 1, any $1 \leq k \leq c(\frac{\log N}{\log\log N})^{1/3}$ and any $0 \leq j \leq k$, we have
$$ \E\Biggl| \int_{0}^{1} \left(\sum_{n \leq N} a_n f(n) e(n\theta) \right)^j \left(\overline{\sum_{n \leq N} a_n f(n) e(n\theta)} \right)^k d\theta - k! (\sum_{n \leq N} |a_n|^2)^k \textbf{1}_{j=k} \Biggr|^2 \ll \frac{N^{j+k}}{N^{1/15k}} , $$
where $\textbf{1}$ denotes the indicator function.

Under the same conditions, and provided the $a_n$ are real, we have
$$ \E\Biggl| \int_{0}^{1} \left(\sum_{n \leq N} a_n f(n) \cos(2\pi n\theta) \right)^k d\theta - \frac{k!}{(k/2)! 2^{k/2}} (\frac{1}{2} \sum_{n \leq N} |a_n|^2)^{k/2} \textbf{1}_{k \; \text{even}} \Biggr|^2 \ll \frac{N^{k}}{N^{1/15k}} , $$
and the same when $\sum_{n \leq N} a_n f(n) \cos(2\pi n\theta)$ is replaced by $\sum_{n \leq N} a_n f(n) \sin(2\pi n\theta)$.
\end{prob2}

\begin{proof}[Proof of Probability Result \ref{rmfimport}]
When $a_n = 1$ for all $n$, the first statement follows immediately from Theorem 1.1 of Benatar, Nishry and Rodgers~\cite{bennishrodg}, after adjusting for the rescaling of the sums by $1/\sqrt{N}$ that they perform but we do not, and handling the easy $j=0$ case that they omit. For general $a_n$, one can check that the proof of Theorem 1.1 transfers over straightforwardly, since the diagonal contribution to the integral (coming when $j=k$, from summands in $\left(\sum_{n \leq N} a_n f(n) e(n\theta) \right)^k$ that are a permutation of the summands in $\left(\overline{\sum_{n \leq N} a_n f(n) e(n\theta)} \right)^k$) has the acceptable form
$$ k! \sum_{n_1, ..., n_k \leq N} |a_{n_1}|^2 ... |a_{n_k}|^2 + O\Biggl( k! \sum_{\substack{n_1, ..., n_k \leq N, \\ n_i \; \text{not all distinct}}} |a_{n_1}|^2 ... |a_{n_k}|^2 \Biggr) = k! (\sum_{n \leq N} |a_n|^2)^k + O(k! k^2 N^{k-1}) , $$
and all off-diagonal contributions continue to satisfy the point-counting bounds of Benatar, Nishry and Rodgers~\cite{bennishrodg} (since the weights $a_n$ are bounded in absolute value by 1).

To deduce the second statement, by writing $\cos(2\pi n\theta) = \frac{e(n\theta) + e(-n\theta)}{2}$ and expanding the $k$-th power we can rewrite $\int_{0}^{1} \left(\sum_{n \leq N} a_n f(n) \cos(2\pi n\theta) \right)^k d\theta$ as a weighted sum of $k+1$ terms of the form $\int_{0}^{1} \left(\sum_{n \leq N} a_n f(n) e(n\theta) \right)^j \left(\overline{\sum_{n \leq N} a_n f(n) e(n\theta)} \right)^{k-j} d\theta$. Here we use the fact that $a_n$ and $f(n)$ are real valued, and so $\sum_{n \leq N} a_n f(n) e(-n\theta) = \overline{\sum_{n \leq N} a_n f(n) e(n\theta)}$. If $k$ is even, then the term with $j = k/2$ is
$$ \frac{1}{2^k} {{k}\choose{k/2}} \int_{0}^{1} \left(\sum_{n \leq N} a_n f(n) e(n\theta) \right)^{k/2} \left(\overline{\sum_{n \leq N} a_n f(n) e(n\theta)} \right)^{k/2} d\theta , $$
and from the first part of Probability Result \ref{rmfimport} we obtain a corresponding ``main term'' $\frac{1}{2^k} {{k}\choose{k/2}} (k/2)! (\sum_{n \leq N} |a_n|^2)^{k/2} = \frac{k!}{(k/2)! 2^{k/2}} (\frac{1}{2} \sum_{n \leq N} |a_n|^2)^{k/2}$, as desired. No other values of $j$ produce any main terms. So using the first part of Probability Result \ref{rmfimport}, along with Minkowski's inequality, we obtain
\begin{eqnarray}
&& \sqrt{ \E\Biggl| \int_{0}^{1} \left(\sum_{n \leq N} a_n f(n) \cos(2\pi n\theta) \right)^k d\theta - \frac{k!}{(k/2)! 2^{k/2}} (\frac{1}{2} \sum_{n \leq N} |a_n|^2)^{k/2} \textbf{1}_{k \; \text{even}} \Biggr|^2 } \nonumber \\
& \leq & \sum_{j=0}^{k} \frac{1}{2^k} {{k}\choose{j}} \sqrt{ \E\Biggl| \int_{0}^{1} \left(\sum_{n \leq N} a_n f(n) e(n\theta) \right)^j \left(\overline{\sum_{n \leq N} a_n f(n) e(n\theta)} \right)^{k-j} - j! (\sum_{n \leq N} |a_n|^2)^j \textbf{1}_{j=k-j} \Biggr|^2 } \nonumber \\
& \ll & \sum_{j=0}^{k} \frac{1}{2^k} {{k}\choose{j}} \sqrt{ \frac{N^{k}}{N^{1/15k}} } = \sqrt{ \frac{N^{k}}{N^{1/15k}} } . \nonumber
\end{eqnarray}
Squaring both sides yields the second part of Probability Result \ref{rmfimport}.

To handle $\sum_{n \leq N} a_n f(n) \sin(2\pi n\theta)$, one proceeds in the same way writing $\sin(2\pi n\theta) = \frac{e(n\theta) - e(-n\theta)}{2i}$, and noting that if $k$ is even then the term with $j = k/2$ again produces a main term $\frac{1}{(2i)^k} {{k}\choose{k/2}} (-1)^{k/2} (k/2)! (\sum_{n \leq N} |a_n|^2)^{k/2} = \frac{1}{2^k} {{k}\choose{k/2}} (k/2)! (\sum_{n \leq N} |a_n|^2)^{k/2} = \frac{k!}{(k/2)! 2^{k/2}} (\frac{1}{2} \sum_{n \leq N} |a_n|^2)^{k/2}$.
\end{proof}

In the Steinhaus case, the estimates we require cannot be read so immediately out of the work of Benatar, Nishry and Rodgers~\cite{bennishrodg}, but we can extract suitable results by adapting their proofs.

\begin{prob3}\label{rmfsteinhausimport}
Let $f(n)$ be a Steinhaus random multiplicative function. Then uniformly for any large $N$, any coefficients $(a_n)_{n \leq N}$ bounded in absolute value by 1, any $1 \leq k \leq c(\frac{\log N}{\log\log N})^{1/3}$ and any $0 \leq j \leq k$, we have
$$ \E\Biggl| \int_{0}^{1} \left(\sum_{n \leq N} a_n f(n) \cos(2\pi n\theta) \right)^j \left(\overline{\sum_{n \leq N} a_n f(n) \cos(2\pi n\theta)} \right)^k d\theta - k! (\frac{1}{2} \sum_{n \leq N} |a_n|^2)^k \textbf{1}_{j=k} \Biggr|^2 \ll \frac{N^{j+k}}{N^{1/15k}} , $$
where $\textbf{1}$ denotes the indicator function.

The same is true when $\sum_{n \leq N} a_n f(n) \cos(2\pi n\theta)$ is replaced by $\sum_{n \leq N} a_n f(n) \sin(2\pi n\theta)$.
\end{prob3}

The proof of Probability Result \ref{rmfsteinhausimport} will rest on the following three Claims.

\begin{cl1}\label{claimsignssquare}
Let $N \in \N$ be large, let $0 \leq j \leq J$ and $1 \leq k \leq K$, and let $\mathcal{A}$ denote the set of all tuples $(m_1, ..., m_J, n_1, ..., n_K) \in \{1, ..., N\}^{J+K}$ satisfying
$$ \sum_{i=1}^{j} m_i - \sum_{i=j+1}^{J} m_i = \sum_{i=1}^{k} n_i - \sum_{i=k+1}^{K} n_i  , \;\;\;\;\; \text{and} \;\;\;\;\; \prod_{i=1}^{J} m_i \cdot \prod_{i=1}^{K} n_i \; \text{is a square} . $$
Given such a tuple, let $\textbf{m}^{w}$ denote the weighted set obtained from $(m_1, ..., m_J)$ by counting each element $m$ with weight $w(m) = \#\{1 \leq i \leq j : m_i = m\} - \#\{j+1 \leq i \leq J : m_i = m\}$ (and discarding any elements whose weight turns out to be zero). Let $\textbf{n}^{w}$ denote the analogous weighted set obtained from $(n_1, ..., n_K)$.

Then the number of tuples in $\mathcal{A}$ that do {\em not} satisfy $\textbf{m}^{w} = \textbf{n}^{w}$ (i.e. equality of the elements of the sets and of their weights) is $\ll N^{(J+K)/2} \exp\{-\frac{\log N}{3(J+K)} + O((J+K)^2 (\log(J+K) + \log\log N))\}$.
\end{cl1}

\begin{proof}[Proof of Claim \ref{claimsignssquare}]
We can rearrange the conditions defining $\mathcal{A}$ into the form
$$ \sum_{i=1}^{j} m_i + \sum_{i=k+1}^{K} n_i = \sum_{i=1}^{k} n_i + \sum_{i=j+1}^{J} m_i , \;\;\;\;\; \text{and} \;\;\;\;\; \prod_{i=1}^{J} m_i \cdot \prod_{i=1}^{K} n_i \; \text{is a square} . $$
And the relation $\textbf{m}^{w} = \textbf{n}^{w}$ is equivalent to saying that the tuple $(m_1, ..., m_j, n_{k+1}, ..., n_K)$ now occurring on the left is a permutation of the tuple $(n_1, ..., n_k, m_{j+1}, ..., m_J)$ on the right. So by Lemma 3.2 of Benatar, Nishry and Rodgers~\cite{bennishrodg} (writing the bound in the slightly more precise form from display (3.10) in their proof), the number of tuples in $\mathcal{A}$ with $\textbf{m}^{w} \neq \textbf{n}^{w}$ is indeed $\ll N^{(J+K)/2} \exp\{-\frac{\log N}{3(J+K)} + O((J+K)^2 (\log(J+K) + \log\log N))\}$.
\end{proof}

If we replace the condition that $\prod_{i=1}^{J} m_i \cdot \prod_{i=1}^{K} n_i$ is a square by the stronger condition that $\prod_{i=1}^{J} m_i = \prod_{i=1}^{K} n_i$, then we can obtain another relationship between $(m_1, ..., m_J)$ and $(n_1, ..., n_K)$ (for all except a small collection of tuples).

\begin{cl2}\label{claimsignsstein}
Let $N \in \N$ be large, let $0 \leq j \leq J$ and $1 \leq k \leq K$, and let $\mathcal{B}$ denote the set of all tuples $(m_1, ..., m_J, n_1, ..., n_K) \in \{1, ..., N\}^{J+K}$ satisfying
$$ \sum_{i=1}^{j} m_i - \sum_{i=j+1}^{J} m_i = \sum_{i=1}^{k} n_i - \sum_{i=k+1}^{K} n_i  , \;\;\;\;\; \text{and} \;\;\;\;\; \prod_{i=1}^{J} m_i = \prod_{i=1}^{K} n_i . $$
Then the number of tuples in $\mathcal{B}$ for which $(m_1, ..., m_J)$ is {\em not} a permutation of $(n_1, ..., n_K)$ is $\ll N^{(J+K)/2} \exp\{-\frac{\log N}{3(J+K)} + O((J+K)^2 (\log(J+K) + \log\log N))\}$.
\end{cl2}

\begin{proof}[Proof of Claim \ref{claimsignsstein}]
Since any tuple in $\mathcal{B}$ also belongs to the set $\mathcal{A}$ from Claim \ref{claimsignssquare}, we may restrict attention to tuples satisfying $\textbf{m}^{w} = \textbf{n}^{w}$. We shall analyse these by investigating the number $s$ of common elements (counted with multiplicity) between the multisets $\{m_1, ..., m_j\}$ and $\{m_{j+1}, ..., m_J\}$, and the number $t$ of common elements between $\{n_1, ..., n_k\}$ and $\{n_{k+1}, ..., n_K\}$. If $s=t=0$ then the relation $\textbf{m}^{w} = \textbf{n}^{w}$ implies that $(m_1, ..., m_J)$ is a permutation of $(n_1, ..., n_K)$, so we may ignore this case and assume that $s+t \geq 1$.

After possibly reordering some of the $m_i$ and $n_i$ (which at worst will multiply our final bounds by an acceptable factor of $(J+K)!$), we may assume that $m_i = m_{j+i}$ for all $1 \leq i \leq s$ and that $n_i = n_{k+i}$ for all $1 \leq i \leq t$. This leaves $J+K-2(s+t)$ other components of $(m_1, ..., m_J, n_1, ..., n_K)$. The relation $\textbf{m}^{w} = \textbf{n}^{w}$ implies these $J+K-2(s+t)$ latter components must consist of $(1/2)(J+K-2(s+t))$ components $m_i$ (for which there are $\leq N^{(J+K)/2 - s - t}$ possibilities), and $(1/2)(J+K-2(s+t))$ components $n_i$ that are a permutation of the $m_i$. Meanwhile, note that
$$ \textbf{m}^{w} = \textbf{n}^{w} \;\;\ \text{and} \;\; \prod_{i=1}^{J} m_i = \prod_{i=1}^{K} n_i \;\;\;\;\; \Rightarrow \;\;\;\;\; \prod_{i=1}^{s} m_i^2 = \prod_{i=1}^{t} n_i^2 \;\;\;\;\; \Rightarrow \;\;\;\;\; \prod_{i=1}^{s} m_i = \prod_{i=1}^{t} n_i . $$
Then standard calculations with iterated divisor functions $d_{\alpha}(\cdot)$ (see e.g. section 3.1 of Benatar, Nishry and Rodgers~\cite{bennishrodg}) show the number of possibilities for $m_1, ..., m_s, n_1, ..., n_t$ is
\begin{eqnarray}
& \leq & \min\{\sum_{m_1, ..., m_s \leq N} d_{t}(m_1 ... m_s) , \sum_{n_1, ..., n_t \leq N} d_{s}(n_1 ... n_t) \} \nonumber \\
& \leq & \min\{(\sum_{m \leq N} d_{t}(m))^s , (\sum_{n \leq N} d_{s}(n) )^t \} \leq N^{\min\{s,t\}} (2\log N)^{st} \leq N^{\min\{s,t\}} (2\log N)^{(J+K)^2} . \nonumber
\end{eqnarray}
So for given $s$ and $t$, our total number of possible tuples is
$$ \leq (J+K)! \cdot N^{(J+K)/2 - s - t} \cdot N^{\min\{s,t\}} (2\log N)^{(J+K)^2} \leq N^{(J+K)/2 - \max\{s,t\}} e^{O((J+K)\log(J+K) + (J+K)^2 \log\log N)} . $$

Summing over all $s+t \geq 1$ gives a more than acceptable final contribution.
\end{proof}

We shall also require a slightly more complicated, ``doubled up'' version of Claim \ref{claimsignsstein}.

\begin{cl3}\label{claimsignsdouble}
Let $N \in \N$ be large, let $0 \leq j \leq J$ and $1 \leq k \leq K$, and let $\mathcal{C}$ denote the set of all tuples $(m_1^{(1)}, ..., m_J^{(1)}, n_1^{(1)}, ..., n_K^{(1)}, m_1^{(2)}, ..., m_J^{(2)}, n_1^{(2)}, ..., n_K^{(2)}) \in \{1, ..., N\}^{2(J+K)}$ satisfying
$$ \sum_{i=1}^{j} m_i^{(l)} - \sum_{i=j+1}^{J} m_i^{(l)} = \sum_{i=1}^{k} n_i^{(l)} - \sum_{i=k+1}^{K} n_i^{(l)} \;\; \forall l \in \{1,2\} , \;\;\;\;\; \text{and} \;\;\;\;\; \prod_{i=1}^{J} m_i^{(1)} \cdot \prod_{i=1}^{K} n_i^{(2)} = \prod_{i=1}^{J} m_i^{(2)} \cdot \prod_{i=1}^{K} n_i^{(1)} . $$
Then the number of tuples in $\mathcal{C}$ for which $(m_1^{(1)}, ..., m_J^{(1)})$ is {\em not} a permutation of $(n_1^{(1)}, ..., n_K^{(1)})$, or $(m_1^{(2)}, ..., m_J^{(2)})$ is {\em not} a permutation of $(n_1^{(2)}, ..., n_K^{(2)})$, is $\ll N^{J+K} \exp\{-\frac{\log N}{6(J+K)} + O((J+K)^2 (\log(J+K) + \log\log N))\}$.
\end{cl3}

\begin{proof}[Proof of Claim \ref{claimsignsdouble}]
With an obvious adaptation of the notation from Claim \ref{claimsignssquare}, we may restrict attention to tuples in $\mathcal{C}$ that satisfy $\textbf{m}^{(1)w} = \textbf{n}^{(1)w}$ and $\textbf{m}^{(2)w} = \textbf{n}^{(2)w}$. For if there is some element whose weight is (say) greater in $\textbf{m}^{(1)}$ than in $\textbf{n}^{(1)}$, and whose weight is at least as great in $\textbf{m}^{(2)}$ as in $\textbf{n}^{(2)}$, then (adding our equations\footnote{Note that in general we might need to swap the roles of $\textbf{m}^{(2)}$ and $\textbf{n}^{(2)}$, in other words add the left hand side of our $l=1$ equation to the right hand side of our $l=2$ equation.} for $l=1,2$) we get
$$ \sum_{i=1}^{j} m_i^{(1)} + \sum_{i=1}^{j} m_i^{(2)} - \sum_{i=j+1}^{J} m_i^{(1)} - \sum_{i=j+1}^{J} m_i^{(2)} = \sum_{i=1}^{k} n_i^{(1)} + \sum_{i=1}^{k} n_i^{(2)} - \sum_{i=k+1}^{K} n_i^{(1)} - \sum_{i=k+1}^{K} n_i^{(2)} , $$
$$ \text{and} \;\;\;\;\; \prod_{i=1}^{J} m_i^{(1)} \cdot \prod_{i=1}^{J} m_i^{(2)} \cdot \prod_{i=1}^{K} n_i^{(1)} \cdot \prod_{i=1}^{K} n_i^{(2)} \; \text{is a square} , $$
where the weighted set corresponding to the concatenated tuple of $m$ terms on the left will be unequal to the weighted set corresponding to the $n$ terms on the right. Thus Claim \ref{claimsignssquare} implies that the number of such ``bad'' tuples is $\ll N^{J+K} \exp\{-\frac{\log N}{6(J+K)} + O((J+K)^2 (\log(J+K) + \log\log N))\}$.

For those ``good'' tuples where $\textbf{m}^{(1)w} = \textbf{n}^{(1)w}$ and $\textbf{m}^{(2)w} = \textbf{n}^{(2)w}$, we may conclude similarly as in the proof of Claim \ref{claimsignsstein}. Thus if $s^{(1)}$ denotes the number of common elements (counted with multiplicity) between $\{m_1^{(1)}, ..., m_j^{(1)}\}$ and $\{m_{j+1}^{(1)}, ..., m_J^{(1)}\}$, similarly for $t^{(1)}, s^{(2)}, t^{(2)}$, then we may ignore the case where all of these are zero, and otherwise our total number of possible tuples is
\begin{eqnarray}
& \ll & (2(J+K))! \cdot N^{(J+K)/2 - s^{(1)} - t^{(1)}} N^{(J+K)/2 - s^{(2)} - t^{(2)}} \cdot N^{\min\{s^{(1)} + t^{(2)}, s^{(2)} + t^{(1)}\}} (2\log N)^{(s^{(1)} + t^{(2)})(s^{(2)} + t^{(1)})} \nonumber \\
& \ll & N^{J+K- \max\{s^{(1)} + t^{(2)}, s^{(2)} + t^{(1)}\}} e^{O((J+K)\log(J+K) + (J+K)^2 \log\log N)} . \nonumber
\end{eqnarray}
Summing this over all $s^{(1)} + t^{(1)} + s^{(2)} + t^{(2)} \geq 1$ gives an acceptable contribution.
\end{proof}

\begin{proof}[Proof of Probability Result \ref{rmfsteinhausimport}]
Writing $\cos(2\pi n\theta) = \frac{e(n\theta) + e(-n\theta)}{2}$, and attempting to mimic the proof of the Steinhaus case of Theorem 1.1 of Benatar, Nishry and Rodgers~\cite{bennishrodg}, one finds that in place of the linear equations $\sum_{i=1}^{j} m_i = \sum_{i=1}^{k} n_i$ that they encounter we must handle the more general situation where some of the $m_i$ and $n_i$ come with negative signs (arising from the $e(-n\theta)$ terms). Using Claims \ref{claimsignsstein} and \ref{claimsignsdouble} in place of Lemma 3.3 and Corollary 3.5 of Benatar, Nishry and Rodgers~\cite{bennishrodg}, one can bound all the ``off-diagonal'' contributions with the same quality bounds as Benatar, Nishry and Rodgers~\cite{bennishrodg}. Thus it only remains to check that the diagonal contribution to the integral in Probability Result \ref{rmfsteinhausimport} (coming when $j=k$, from summands in $\left(\sum_{n \leq N} a_n f(n) \cos(2\pi n\theta) \right)^k$ that are a permutation of the summands in $\left(\overline{\sum_{n \leq N} a_n f(n) \cos(2\pi n\theta)} \right)^k$) is acceptably close to $k! (\frac{1}{2} \sum_{n \leq N} |a_n|^2)^k$.

But we can write that diagonal contribution as
$$ k! \sum_{n_1, ..., n_k \leq N} |a_{n_1}|^2 ... |a_{n_k}|^2 \int_{0}^{1} \cos^{2}(2\pi n_1 \theta) ... \cos^{2}(2\pi n_k \theta) d\theta + O\Biggl( k! \sum_{\substack{n_1, ..., n_k \leq N, \\ n_i \; \text{not all distinct}}} |a_{n_1}|^2 ... |a_{n_k}|^2 \Biggr) . $$
Since $\cos^{2}(2\pi n\theta) = \frac{1}{2} + \frac{e(2n\theta)}{4} + \frac{e(-2n\theta)}{4}$, we have $\int_{0}^{1} \cos^{2}(2\pi n_1 \theta) ... \cos^{2}(2\pi n_k \theta) d\theta = 1/2^k$ (coming from the term $1/2$ in the expansion of all the factors $\cos^{2}(2\pi n_j \theta)$) except for tuples $n_1, ..., n_k$ satisfying additional linear relations (producing additional contributions from a product of terms $\frac{e(2n_j \theta)}{4} , \frac{e(-2n_j \theta)}{4}$). The total of all such additional contributions, together with the ``big Oh'' term $O\Biggl( k! \sum_{\substack{n_1, ..., n_k \leq N, \\ n_i \; \text{not all distinct}}} |a_{n_1}|^2 ... |a_{n_k}|^2 \Biggr)$, is $\ll k! k^2 N^{k-1} $
\end{proof}

\vspace{12pt}
In order to bring Probability Results \ref{rmfimport} and \ref{rmfsteinhausimport} to bear, we must show that averages of Dirichlet characters behave similarly to averages of random multiplicative functions. When averaging over all characters mod $q$, this will be straightforward (provided we keep sufficient control on the lengths of the sums being averaged) thanks to orthogonality of characters. When averaging only over Legendre symbols $\left(\frac{\cdot}{q}\right)$ with $q$ prime, matters are more subtle, and connected with the distribution of zeros of Dirichlet $L$-functions. Nevertheless there are various approaches that can be applied, for example using the explicit formula for character sums over primes along with results of Siegel and Linnik type on exceptional zeros and (log-free) zero density. Since we will arrange our arguments so that only upper bounds (rather than asymptotic equalities) for character averages are needed, we instead proceed in a different way using the sieve, which will allow quantitatively stronger conclusions about the density of $\mathcal{P}_{H}$ in Theorem 3.

\begin{numb3}[See Lemma 9 of Montgomery and Vaughan~\cite{mvcharmean}, 1979]\label{leglikermf}
Let $f(n)$ be an extended Rademacher random multiplicative function. Then uniformly for any large $Q$, any $N \leq Q$, and any complex coefficients $(\alpha_n)_{n \leq N}$, we have
\begin{eqnarray}
\frac{\log Q}{Q} \sum_{\substack{Q \leq q \leq 2Q, \\ q \; \text{prime}}} \left| \sum_{n \leq N} \alpha_n \left(\frac{n}{q}\right) \right|^2 & \ll & \E\left| \sum_{n \leq N} \alpha_n f(n) \right|^2 + \frac{1}{Q^{0.99}} \Biggl( \sum_{n \leq N} |\alpha_n| \sqrt{s(n)} \Biggr)^2 \nonumber \\
& \ll & \E\left| \sum_{n \leq N} \alpha_n f(n) \right|^2 + \frac{N}{Q^{0.99}} \Biggl( \sum_{n \leq N} |\alpha_n| \Biggr)^2 , \nonumber
\end{eqnarray}
say, where $s(n)$ denotes the squarefree part of $n$ (i.e. $n$ divided by its largest square factor).
\end{numb3}

\begin{proof}[Proof of Number Theory Result \ref{leglikermf}]
This result is very close to Lemma 9 of Montgomery and Vaughan~\cite{mvcharmean}, and would follow by tweaking the argument in Lemmas 4--9 of their paper. For convenience, and since it is neat and fairly short, we outline a self-contained proof here. 

Note first that if $q$ is prime and $n \leq q$, then $\left(\frac{n}{q}\right) = \left(\frac{s(n)}{q}\right)$ where $s(n)$ denotes the squarefree part of $n$. So we can rewrite $\sum_{n \leq N} \alpha_n \left(\frac{n}{q}\right) = \sum_{\substack{s \leq N, \\ s \; \text{squarefree}}} \left(\frac{s}{q}\right) \sum_{\substack{n \leq N, \\ s(n) = s}} \alpha_n$. We also always have $f(n) = f(s(n))$, so it is easy to see that $\E\left| \sum_{n \leq N} \alpha_n f(n) \right|^2 = \E\left| \sum_{\substack{s \leq N, \\ s \; \text{squarefree}}} f(s) \sum_{\substack{n \leq N, \\ s(n) = s}} \alpha_n \right|^2 = \sum_{\substack{s \leq N, \\ s \; \text{squarefree}}} \left| \sum_{\substack{n \leq N, \\ s(n) = s}} \alpha_n \right|^2$.

To execute the proof, the only (possibly) non-obvious step is the introduction of upper bound sieve weights in place of the sum over primes $q$. At the level of precision we are seeking we have much flexibility in our choice of sieve. For example (following the notation of section 3.2 of Montgomery and Vaughan's book~\cite{mv} with $z = Q^{0.005}$ and $P = \prod_{\text{primes} \; p \leq Q^{0.005}} p$), we can use Selberg sieve weights $\lambda_e = \lambda_{e}^{+}$ satisfying $\lambda_e = 0$ whenever $e > Q^{0.01}$, and $\sum_{e} |\lambda_e| \ll \frac{Q^{0.01}}{\log^{2}Q}$, and $\sum_{e|q} \lambda_e \geq \textbf{1}_{p|q \Rightarrow p > Q^{0.005}}$ for all $q$, and $\sum_{Q \leq q \leq 2Q} \sum_{e|q} \lambda_e \ll \frac{Q}{\log Q}$. Thus we have
\begin{eqnarray}
\frac{\log Q}{Q} \sum_{\substack{Q \leq q \leq 2Q, \\ q \; \text{prime}}} \left| \sum_{n \leq N} \alpha_n \left(\frac{n}{q}\right) \right|^2 & \leq & \frac{\log Q}{Q} \sum_{\substack{Q \leq q \leq 2Q, \\ q \; \text{odd}}} (\sum_{e|q} \lambda_e ) \Biggl| \sum_{\substack{s \leq N, \\ s \; \text{squarefree}}} \left(\frac{s}{q}\right) \sum_{\substack{n \leq N, \\ s(n) = s}} \alpha_n \Biggr|^2 \nonumber \\
& = & \sum_{\substack{s_1 , s_2 \leq N, \\ \text{squarefree}}} (\sum_{\substack{n \leq N, \\ s(n) = s_1}} \alpha_n) \overline{(\sum_{\substack{n \leq N, \\ s(n) = s_2}} \alpha_n)} \frac{\log Q}{Q} \sum_{\substack{Q \leq q \leq 2Q, \\ q \; \text{odd}}} (\sum_{e|q} \lambda_e ) \left(\frac{s_1 s_2}{q}\right) . \nonumber
\end{eqnarray}
Here $\left(\frac{s}{q}\right)$ should be understood to mean the Jacobi symbol, which is well defined for all odd $q$ and all $s$, and agrees with the Legendre symbol when $q$ is prime. See section 9.3 of Montgomery and Vaughan's book~\cite{mv}, for example.

The contribution from the diagonal summands $s_1 = s_2$ is
\begin{eqnarray}
\sum_{\substack{s \leq N, \\ \text{squarefree}}} |\sum_{\substack{n \leq N, \\ s(n) = s}} \alpha_n |^2 \frac{\log Q}{Q} \sum_{\substack{Q \leq q \leq 2Q, \\ q \; \text{odd}}} (\sum_{e|q} \lambda_e ) \left(\frac{s^2}{q}\right) & \leq & \sum_{\substack{s \leq N, \\ \text{squarefree}}} |\sum_{\substack{n \leq N, \\ s(n) = s}} \alpha_n |^2 \frac{\log Q}{Q} \sum_{Q \leq q \leq 2Q} (\sum_{e|q} \lambda_e ) \nonumber \\
& \ll & \sum_{\substack{s \leq N, \\ \text{squarefree}}} |\sum_{\substack{n \leq N, \\ s(n) = s}} \alpha_n |^2 = \E\left| \sum_{n \leq N} \alpha_n f(n) \right|^2 , \nonumber
\end{eqnarray}
which is acceptable.

If $s_1 \neq s_2$ are squarefree, then $s_1 s_2$ is not a perfect square, and so the mapping $q \mapsto \left(\frac{s_1 s_2}{q}\right)$ is a non-principal Dirichlet character of conductor at most $4s_1 s_2$. Hence we can bound the contribution from $s_1 \neq s_2$ by
\begin{eqnarray}
&& \frac{\log Q}{Q} \sum_{\substack{e \leq Q^{0.01}, \\ e \; \text{odd}}} |\lambda_e| \sum_{\substack{s_1 \neq s_2 \leq N, \\ \text{squarefree}}} |\sum_{\substack{n \leq N, \\ s(n) = s_1}} \alpha_n| |\sum_{\substack{n \leq N, \\ s(n) = s_2}} \alpha_n | \Biggl|\sum_{\substack{Q \leq q \leq 2Q, \\ q \; \text{odd}, \\ e|q}} \left(\frac{s_1 s_2}{q}\right) \Biggr| \nonumber \\
& \ll & \frac{\log Q}{Q} \sum_{\substack{e \leq Q^{0.01} , \\ e \; \text{odd}}} |\lambda_e| \sum_{\substack{s_1 \neq s_2 \leq N, \\ \text{squarefree}}} |\sum_{\substack{n \leq N, \\ s(n) = s_1}} \alpha_n| |\sum_{\substack{n \leq N, \\ s(n) = s_2}} \alpha_n | \sqrt{s_1 s_2} \log(s_1 s_2) , \nonumber
\end{eqnarray}
where the second line follows using the P\'olya--Vinogradov inequality (see e.g. section 9.4 of Montgomery and Vaughan~\cite{mv}) and a little manipulation. (Note that because we switched to sums with sieve weights rather than sums over primes, here we finally obtained character sums over (essentially) all integers $q$ in an interval, for which we have the strong P\'olya--Vinogradov bound.) Since our weights $\lambda_e$ satisfy $\sum_{e} |\lambda_e| \ll \frac{Q^{0.01}}{\log^{2}Q}$, one can check that this expression is also acceptably small.
\end{proof}

% SECTION 6 %%%%%%%%%%%%%%%%%%%%%%%%%%%%%%
\section{Proof of Theorem 3}
In this section we shall prove Theorem 3, our positive ``almost all'' result for real characters. The proof splits into four parts: firstly we shall reduce the problem to one about the distribution of sufficiently short exponential sums (this part will also be applicable when handling the complex case in Theorem 4); secondly we show that it will suffice to bound mean square averages (over characters) of moment related objects involving those exponential sums; thirdly, we perform a technical ``netting'' step allowing us to treat $q/H(q)$ as constant on dyadic ranges $Q \leq q \leq 2Q$, so that we can perform the desired averages over $q$ (this is only needed in the real case); and finally we complete the analysis using Probability Result \ref{rmfimport} and Number Theory Result \ref{leglikermf}. 

\subsection{Reduction to short partial Fourier series}\label{fourierred}
If $\chi$ is an even non-principal Dirichlet character mod $q$ (so that $\chi(-1)=1$, and therefore $\chi(-k) = \chi(k)$ for all $k$), then we can rewrite the P\'{o}lya Fourier expansion \eqref{polyaq} in the form
\begin{eqnarray}
\frac{S_{\chi, H}(X)}{\sqrt{H}} & = & \frac{\tau(\chi)}{2\pi i \sqrt{H}} \sum_{1 \leq k < q/2} \frac{\overline{\chi}(k)}{k} \left(e(\frac{kX}{q}) (e(kH/q) - 1) - e(-\frac{kX}{q}) (e(-kH/q) - 1) \right) + O(\frac{\log q}{\sqrt{H}}) \nonumber \\
& = & \frac{\tau(\chi)}{\pi \sqrt{H}} \sum_{1 \leq k < q/2} \frac{\overline{\chi}(k)}{k} \left(\sin(2\pi k(X+H)/q) - \sin(2\pi k X/q) \right) + O(\frac{\log q}{\sqrt{H}}) \nonumber \\
& = & \frac{2\tau(\chi)}{\pi \sqrt{H}} \sum_{1 \leq k < q/2} \frac{\overline{\chi}(k) \sin(\pi k H/q)}{k} \cos(\pi k (2X+H)/q) + O(\frac{\log q}{\sqrt{H}}) . \nonumber
\end{eqnarray}

The sum over $k$ here would be too long for our subsequent calculations, in particular to allow the computation of its high moments. However, if $1 \leq k_1 , k_2 < q/2$ and if $X \in \{0,1,...,q-1\}$ is uniformly random then we have $\E \cos(\pi k_1 (2X+H)/q) \cos(\pi k_2 (2X+H)/q) = (1/2) \textbf{1}_{k_1 = k_2}$, and so
\begin{eqnarray}
&& \E\left| \frac{2\tau(\chi)}{\pi \sqrt{H}} \sum_{(q/H)\log(q/H) \leq k < q/2} \frac{\overline{\chi}(k) \sin(\pi k H/q)}{k} \cos(\pi k (2X+H)/q) \right|^2 \nonumber \\
& = & \frac{2q}{\pi^2 H} \sum_{(q/H)\log(q/H) \leq k < q/2} \frac{\sin^{2}(\pi k H/q)}{k^2} \ll \frac{1}{\log(q/H)} , \nonumber
\end{eqnarray}
which tends to zero as $q \rightarrow \infty$ under the conditions of Theorem 3 (or Theorem 4). It follows that the part of the sum with $k \geq (q/H)\log(q/H)$ tends to zero in probability, for {\em any} choice of $\chi$, so may be ignored in our investigation of the limiting distribution.

The form of the function $\cos(\pi k (2X+H)/q)$, with $X \in \{0,1,...,q-1\}$ uniformly random, is a bit ungainly. However, under the conditions of Theorems 3 and 4 it turns out we can replace this by $\cos(2 \pi k \theta)$, where $\theta \in [0,1]$ is uniformly random. Indeed, if $X \in \{0,1,...,q-1\}$ then for any $\theta \in [\frac{X + H/2}{q} - \frac{1}{2q}, \frac{X + H/2}{q} + \frac{1}{2q}]$ mod 1 we get
\begin{eqnarray}
&& \left| \frac{2\tau(\chi)}{\pi \sqrt{H}} \sum_{1 \leq k < (q/H)\log(q/H)} \frac{\overline{\chi}(k) \sin(\pi k H/q)}{k} (\cos(\pi k (2X + H)/q) - \cos(2\pi k\theta)) \right| \nonumber \\
& \ll & \frac{\sqrt{q}}{\sqrt{H}} \sum_{1 \leq k < (q/H)\log(q/H)} \frac{|\sin(\pi k H/q)|}{k} \cdot \frac{k}{q} \ll \frac{(q/H)^{3/2} \log(q/H)}{q} , \nonumber
\end{eqnarray}
which tends to zero (deterministically) as $q \rightarrow \infty$. Here we mildly use our assumption that $\frac{\log(q/H)}{\log q} \rightarrow 0$ as $q \rightarrow \infty$. Since choosing $X \in \{0,1,...,q-1\}$ uniformly at random, and then choosing $\theta \in [\frac{X + H/2}{q} - \frac{1}{2q}, \frac{X + H/2}{q} + \frac{1}{2q}]$ mod 1 uniformly at random, is exactly the same thing as choosing $\theta \in [0,1]$ uniformly at random, we only need to consider the latter process.

In summary, for even characters $\chi$ it will suffice to prove Theorem 3 (and Theorem 4) with $S_{\chi, H}(X)/\sqrt{H}$ replaced by
$$ \frac{2\tau(\chi)}{\pi \sqrt{H}} \sum_{1 \leq k < (q/H)\log(q/H)} \frac{\overline{\chi}(k) \sin(\pi k H/q)}{k} \cos(2\pi k\theta) , \;\;\;\;\; \theta \sim \text{Uni}[0,1] . $$
For odd characters $\chi$, where $\chi(-1)=-1$ and therefore $\chi(-k) = -\chi(k)$ for all $k$, one similarly ends up with $\frac{2 \tau(\chi)}{\pi i \sqrt{H}} \sum_{1 \leq k < (q/H)\log(q/H)} \frac{\overline{\chi}(k) \sin(\pi k H/q)}{k} \sin(2\pi k\theta)$. The treatment of either sum will be exactly similar, so for simplicity we shall focus on the cosine case. Note that we need not distinguish between even and odd characters in our subsequent calculations, because if $\frac{2\tau(\chi)}{\pi \sqrt{H}} \sum_{1 \leq k < (q/H)\log(q/H)} \frac{\overline{\chi}(k) \sin(\pi k H/q)}{k} \cos(2\pi k\theta)$ has the desired Gaussian limiting distribution for ``almost all'' choices of $\chi$ (in the sense of Theorems 3 and 4) then, in particular, it has the desired limiting distribution for almost all choices of even $\chi$, similarly for $\frac{2 \tau(\chi)}{\pi i \sqrt{H}} \sum_{1 \leq k < (q/H)\log(q/H)} \frac{\overline{\chi}(k) \sin(\pi k H/q)}{k} \sin(2\pi k\theta)$.

\vspace{12pt}
Note also that when $\chi = \left(\frac{\cdot}{q}\right)$ is real, one has $\tau(\chi) = \sqrt{q}$ if $\chi$ is even (which occurs when $q \equiv 1$ mod 4), and one has $\tau(\chi) = i\sqrt{q}$ if $\chi$ is odd (which occurs when $q \equiv 3$ mod 4). See chapter 9.3 of Montgomery and Vaughan~\cite{mv}. Inserting these expressions above, we see that when proving Theorem 3 we can work with
\begin{equation}\label{polishedsum}
\frac{2\sqrt{q}}{\pi \sqrt{H}} \sum_{1 \leq k < (q/H)\log(q/H)} \frac{\overline{\chi}(k) \sin(\pi k H/q)}{k} \cos(2\pi k\theta) , \;\;\;\;\; \theta \sim \text{Uni}[0,1] 
\end{equation}
and with $\frac{2 \sqrt{q}}{\pi \sqrt{H}} \sum_{1 \leq k < (q/H)\log(q/H)} \frac{\overline{\chi}(k) \sin(\pi k H/q)}{k} \sin(2\pi k\theta)$. These sums are visibly real-valued when $\chi$ is real. And in fact we are free to work with these sums when proving Theorem 4 as well, where we know that $|\tau(\chi)| = \sqrt{q}$, but it is harder to say a lot about the argument of $\tau(\chi)$. That is because in Theorem 4 the target distribution $Z_1 + iZ_2$ is rotationally invariant, so if this is the limiting distribution of e.g. \eqref{polishedsum} for ``almost all'' choices of $\chi$ then it remains the limiting distribution of $\frac{2\tau(\chi)}{\pi \sqrt{H}} \sum_{1 \leq k < (q/H)\log(q/H)} \frac{\overline{\chi}(k) \sin(\pi k H/q)}{k} \cos(2\pi k\theta)$ for the same $\chi$.

\subsection{Working with moments}\label{momentsmanipsec}
The method of moments for proving distributional convergence is discussed in a general context in e.g. chapter 5.8.4 of Gut~\cite{gut}. In particular, our $N(0,1)$ target distribution is determined by its moments, which are $\E N(0,1)^j = (1/\sqrt{2\pi}) \int_{-\infty}^{\infty} z^j e^{-z^{2}/2} dz = \textbf{1}_{j \; \text{even}} \frac{j!}{2^{j/2} (j/2)!}$ where $\textbf{1}$ denotes the indicator function. So in view of \eqref{polishedsum}, to prove Theorem 3 it would suffice to show that there exists a subsequence $\mathcal{P}_{H}$ of primes, with the density claimed in the theorem, such that for all $j \in \N$ we have
$$ \int_{0}^{1} \left( \frac{2\sqrt{q}}{\pi \sqrt{H}} \sum_{1 \leq k < (q/H)\log(q/H)} \frac{\left(\frac{k}{q}\right) \sin(\pi k H/q)}{k} \cos(2\pi k\theta) \right)^{j} d\theta \rightarrow \frac{j! \textbf{1}_{j \; \text{even}}}{(j/2)! 2^{j/2}} \;\;\;\;\; \text{as} \; q \rightarrow \infty , \; q \in \mathcal{P}_{H} . $$
(Actually we must also prove this for $\frac{2 \sqrt{q}}{\pi \sqrt{H}} \sum_{1 \leq k < (q/H)\log(q/H)} \frac{\left(\frac{k}{q}\right) \sin(\pi k H/q)}{k} \sin(2\pi k\theta)$, but this will be exactly similar to the cosine case, so we shall only discuss the latter.)

Rewriting slightly, if we set $a_k = a_{k,q,H} := \frac{q \sin(\pi k H/q)}{\pi H k}$ then we want to show the existence of $\mathcal{P}_{H}$ such that, for each fixed $j \in \N$, we have
$$ \left(\frac{4H}{q} \right)^{j/2} \left| \int_{0}^{1} \left( \sum_{k < (q/H)\log(q/H)} a_k \left(\frac{k}{q}\right) \cos(2\pi k\theta) \right)^{j} - \frac{j! \textbf{1}_{j \; \text{even}}}{(j/2)! 2^{j/2}} (\frac{q}{4H})^{j/2} \right| \rightarrow 0 \;\;\; \text{as} \; q \rightarrow \infty , \; q \in \mathcal{P}_{H} . $$
Here the coefficients $a_k$ are real, bounded in absolute value by 1 (thanks to the estimate $|\sin x| \leq |x|$), and satisfy
$$ \frac{1}{2} \sum_{k < (q/H)\log(q/H)} |a_k|^2 = \frac{q^2}{2\pi^2 H^2} \sum_{1 \leq k < (q/H)\log(q/H)} \frac{\sin^2(\frac{\pi k H}{q})}{k^2} = \frac{q^2}{4\pi^2 H^2} \sum_{1 \leq |k| < (q/H)\log(q/H)} \frac{\sin^2(\frac{\pi k H}{q})}{k^2} , $$
since $\frac{\sin^2(\pi k H/q)}{k^2}$ is an even function of $k$. Using the fact that $|\sin(\pi k H/q)| = (1/2)|e(k H/2q) - e(- k H/2q)| = (1/2)|e(k H/q) - 1|$, we can rewrite this further as
\begin{eqnarray}
\frac{q^2}{16 \pi^2 H^2} \sum_{1 \leq |k| < (q/H)\log(q/H)} \frac{|e(\frac{k H}{q}) - 1|^2}{k^2} & = & \frac{q^2}{16 \pi^2 H^2} \Biggl(\sum_{1 \leq |k| < q/2} \frac{|e(\frac{k H}{q}) - 1|^2}{k^2} + O(\frac{1}{(q/H)\log(q/H)}) \Biggr) \nonumber \\
& = & \frac{q}{4H} \Biggl( 1 + O(\frac{1}{\log(q/H)}) \Biggr) , \nonumber
\end{eqnarray}
where the final equality uses the calculation of $\frac{q}{(2\pi)^2} \sum_{0 < |k| < q/2} \frac{1}{k^2} |e(kH/q) - 1|^2$ that we performed in section \ref{negproofsec}.

\vspace{12pt}
To finish the proof, for each $Q = 2^r, r \in \N$ we would like to show that the averages
$$ \frac{\log Q}{Q} \sum_{\substack{Q \leq q \leq 2Q, \\ q \; \text{prime}}} \left(\frac{4H}{q} \right)^{j} \left| \int_{0}^{1} \Biggl( \sum_{k < \frac{q\log(q/H)}{H}} a_k \left(\frac{k}{q}\right) \cos(2\pi k\theta) \Biggr)^{j} - \frac{j! \textbf{1}_{j \; \text{even}}}{(j/2)! 2^{j/2}} \Biggl(\frac{1}{2} \sum_{k < \frac{q\log(q/H)}{H}} |a_k|^2 \Biggr)^{j/2} \right|^2 $$
are ``small'', implying that the number of ``bad'' primes $Q \leq q \leq 2Q$ where the summand is large is also small. Unfortunately, as $Q \leq q \leq 2Q$ varies it is not only the character $\chi(k) = \left(\frac{k}{q}\right)$ (which we expect to behave like an extended Rademacher random multiplicative function $f(k)$) that varies here, but also many other terms like $\sqrt{q}/\sqrt{H}$, $\sin(\pi k H/q)$, and the length of the sum over $k$. In other words, the coefficients $a_k = a_{k,q,H}$ may depend a priori on $q/H(q)$ as well as on $k$. Notice this issue will not arise in the non-real case of Theorem 4, where we can average over all characters $\chi$ mod $q$ whilst holding $q$, and therefore $H(q)$ and everything else, fixed. 

\subsection{Controlling the behaviour of $q/H$}\label{fixqH}
To get around the problem just discussed, we will apply a ``netting'' argument to the given function $H(q)$. Given $Q = 2^r, r \in \N$, let us define the small quantity $\eta = \eta_{H,Q} > 0$ by $\eta := \frac{1}{\min_{Q \leq q \leq 2Q} \log(q/H(q))}$, say, and then define a family of functions $H_{n} : [Q,2Q] \rightarrow \R$ in the following way:
$$ H_{n}(q) := \frac{q}{e^{1/\eta + n\eta}} , \;\;\;\;\; 0 \leq n \leq \frac{1}{\eta} (\max_{Q \leq q \leq 2Q} \log(q/H(q)) - \frac{1}{\eta}) . $$

For each $Q \leq q \leq 2Q$, there exists some $n = n(q)$ for which $q/H_{n}(q) \leq q/H(q) \leq e^{\eta} q/H_{n}(q)$, and then
\begin{eqnarray}
&& \int_{0}^{1} \Biggl| \frac{2\sqrt{q}}{\pi \sqrt{H}} \sum_{1 \leq k < (q/H)\log(q/H)} \frac{\left(\frac{k}{q}\right) \sin(\pi k H/q)}{k} \cos(2\pi k\theta) - \nonumber \\
&& - \frac{2\sqrt{q}}{\pi \sqrt{H_n}} \sum_{1 \leq k < (q/H_n)\log(q/H_n)} \frac{\left(\frac{k}{q}\right) \sin(\pi k H_n /q)}{k} \cos(2\pi k\theta) \Biggr|^2 d\theta \nonumber \\
& = & \frac{2}{\pi^2} \sum_{1 \leq k < (q/H_n)\log(q/H_n)} \frac{1}{k^2} |\frac{\sqrt{q}}{\sqrt{H}} \sin(\pi k H /q) - \frac{\sqrt{q}}{\sqrt{H_n}} \sin(\pi k H_n /q)|^2 + O(\frac{1}{\log(q/H)}) \nonumber \\
& \ll & \sum_{1 \leq k < (q/H_n)\log(q/H_n)} \frac{1}{k^2} \min\left\{\frac{(\eta k)^2}{(q/H)}, q/H \right\} + \frac{1}{\log(q/H)} \ll \eta , \nonumber
\end{eqnarray}
where the final line uses the fact that $\frac{d}{dt} \sqrt{t} \sin(\pi k/t) \ll \frac{k}{t^{3/2}}$ for large $t$. The assumptions of Theorem 3 imply that $\eta \rightarrow 0$ as $Q \rightarrow \infty$, and so the difference between \eqref{polishedsum} and the analogous sum involving $H_n$ tends to zero in probability, uniformly for $Q \leq q \leq 2Q$.

Consequently, when proving Theorem 3 it will suffice to work with the particular functions $H_n(q)$, which have the property that $q/H_{n}(q) = Q/H_{n}(Q)$ is constant for all $Q \leq q \leq 2Q$. More precisely: it will {\em suffice to prove} that for each $0 \leq n \leq \frac{1}{\eta} (\max_{Q \leq q \leq 2Q} \log(q/H(q)) - \frac{1}{\eta})$ and all $1 \leq j \leq \min\{\frac{\log^{1/4}(Q/H_{n}(Q))}{50}, \frac{\log Q}{4\log(Q/H_{n}(Q))}\}$, say, we have
\begin{eqnarray}\label{neededafterreg}
&& \frac{\log Q}{Q} \sum_{\substack{Q \leq q \leq 2Q, \\ q \; \text{prime}}} \left| \int_{0}^{1} \Biggl( \sum_{k < \frac{q\log(q/H_n)}{H_n}} a_k \left(\frac{k}{q}\right) \cos(2\pi k\theta) \Biggr)^{j} - \frac{j! \textbf{1}_{j \; \text{even}}}{(j/2)! 2^{j/2}} \Biggl(\frac{1}{2} \sum_{k < \frac{q\log(q/H_n)}{H_n}} |a_k|^2 \Biggr)^{j/2} \right|^2 \nonumber \\
& \ll & \left(\frac{Q}{4H_{n}(Q)} \right)^{j} e^{-2\log^{3/4}(Q/H_{n}(Q))} .
\end{eqnarray}
Note that the coefficients $a_k$ here also depend on $n$, via the value of $q/H_{n}(q)$.

For if \eqref{neededafterreg} holds, then the proportion of ``bad'' primes $Q \leq q \leq 2Q$ for which
\begin{eqnarray}
&& \left(\frac{4H_{n}(q)}{q} \right)^{j} \left| \int_{0}^{1} \Biggl( \sum_{k < \frac{q\log(q/H_n)}{H_n}} a_k \left(\frac{k}{q}\right) \cos(2\pi k\theta) \Biggr)^{j} - \frac{j! \textbf{1}_{j \; \text{even}}}{(j/2)! 2^{j/2}} \Biggl(\frac{1}{2} \sum_{k < \frac{q\log(q/H_n)}{H_n}} |a_k|^2 \Biggr)^{j/2} \right|^2 \nonumber \\
& \gg & e^{-0.1\log^{3/4}(Q/H_{n}(Q))} \nonumber
\end{eqnarray}
must be $\ll e^{-1.9\log^{3/4}(Q/H_{n}(Q))}$, and so the proportion for which this holds for some $1 \leq j \leq \min\{\frac{\log^{1/4}(Q/H_{n}(Q))}{50}, \frac{\log Q}{4\log(Q/H_{n}(Q))}\}$ must be $\ll \log^{1/4}(Q/H_{n}(Q)) e^{-1.9\log^{3/4}(Q/H_{n}(Q))} \ll e^{-1.8\log^{3/4}(Q/H_{n}(Q))}$. Finally, the proportion of primes $Q \leq q \leq 2Q$ that are ``bad'' for some $n$ will be
$$ \ll \sum_{n} e^{-1.8\log^{3/4}(Q/H_{n}(Q))} = \sum_{n} e^{-1.8 (1/\eta + n\eta)^{3/4}} \ll \frac{1}{\eta^2} e^{-1.8 (1/\eta)^{3/4}} \ll e^{- (1/\eta)^{3/4}} = e^{- \min_{Q \leq q \leq 2Q}\log^{3/4}(q/H)} . $$
So if we define our subsequence $\mathcal{P}_H$ of primes by discarding all the bad primes (in the above sense\footnote{Again, to be completely correct we must also discard those primes that will be bad for the corresponding sine series $\sum_{k < \frac{q\log(q/H_n)}{H_n}} a_k \left(\frac{k}{q}\right) \sin(2\pi k\theta)$.}) in each dyadic interval $[Q,2Q]$, then $\mathcal{P}_H$ has the density required in Theorem 3. And since the assumptions of Theorem 3 imply that
$$ \min\{\frac{\log^{1/4}(Q/H_{n}(Q))}{50}, \frac{\log Q}{4\log(Q/H_{n}(Q))}\} \gg \min_{Q \leq q \leq 2Q}\{\log^{1/4}(q/H(q)), \frac{\log q}{\log(q/H(q))}\} \rightarrow \infty $$
as $Q \rightarrow \infty$, the range of $j$ for which the $j$-th moment of $\Biggl( \sum_{k < \frac{q\log(q/H_n)}{H_n}} a_k \left(\frac{k}{q}\right) \cos(2\pi k\theta) \Biggr)^{j}$ approaches the desired Gaussian moment will grow to infinity as $q \rightarrow \infty$ with $q \in \mathcal{P}_H$.

\subsection{The punchline}
It now remains to establish \eqref{neededafterreg}. Recall that for the functions $H_{n}(q)$ from section \ref{fixqH}, the quantity $q/H_{n}(q)$ is constant (depending on $n$) for all $Q \leq q \leq 2Q$, and the coefficients $a_k = \frac{q \sin(\pi k H_{n}/q)}{\pi H_{n} k}$ depend only on $k$ and $n$. Furthermore, by expanding the integral
$$ \int_{0}^{1} \Biggl( \sum_{k < \frac{q\log(q/H_n)}{H_n}} a_k \left(\frac{k}{q}\right) \cos(2\pi k\theta) \Biggr)^{j} = \sum_{m \leq (\frac{q\log(q/H_n)}{H_n})^j} \left(\frac{m}{q}\right) \sum_{\substack{k_1, ..., k_j < \frac{q\log(q/H_n)}{H_n}, \\ k_1 \cdot ... \cdot k_j = m}} \int_{0}^{1} \prod_{i=1}^{j} a_{k_i} \cos(2\pi k_{i} \theta) $$
we see the left hand side in \eqref{neededafterreg} is of the form treated in Number Theory Result \ref{leglikermf}. (The subtracted term $\frac{j! \textbf{1}_{j \; \text{even}}}{(j/2)! 2^{j/2}} \left(\frac{1}{2} \sum_{k < \frac{q\log(q/H_n)}{H_n}} |a_k|^2 \right)^{j/2}$ in \eqref{neededafterreg} may be thought of as part of the coefficient $\alpha_1$ of the trivial Legendre symbol $\left(\frac{1}{q}\right)$.) Note also that provided $Q$ is large enough, we have
$$ (\frac{q\log(q/H_n)}{H_n})^j \leq (\frac{q\log(q/H_n)}{H_n})^{\frac{\log Q}{4\log(q/H_{n})}} \leq Q^{0.26} $$
on the range of $j$ required in \eqref{neededafterreg}. So we may apply Number Theory Result \ref{leglikermf}, and deduce (with a little manipulation of the error term) that
\begin{eqnarray}
&& \frac{\log Q}{Q} \sum_{\substack{Q \leq q \leq 2Q, \\ q \; \text{prime}}} \left| \int_{0}^{1} \Biggl( \sum_{k < \frac{q\log(q/H_n)}{H_n}} a_k \left(\frac{k}{q}\right) \cos(2\pi k\theta) \Biggr)^{j} - \frac{j! \textbf{1}_{j \; \text{even}}}{(j/2)! 2^{j/2}} \Biggl(\frac{1}{2} \sum_{k < \frac{q\log(q/H_n)}{H_n}} |a_k|^2 \Biggr)^{j/2} \right|^2 \nonumber \\
& \ll & \E \left| \int_{0}^{1} \Biggl( \sum_{k < \frac{q\log(q/H_n)}{H_n}} a_k f(k) \cos(2\pi k\theta) \Biggr)^{j} - \frac{j! \textbf{1}_{j \; \text{even}}}{(j/2)! 2^{j/2}} \Biggl(\frac{1}{2} \sum_{k < \frac{q\log(q/H_n)}{H_n}} |a_k|^2 \Biggr)^{j/2} \right|^2 \nonumber \\
&& + \frac{Q^{0.26}}{Q^{0.99}} \Biggl( \Biggl( \sum_{k < \frac{q\log(q/H_n)}{H_n}} |a_k| \Biggr)^{j} + \frac{j! \textbf{1}_{j \; \text{even}}}{(j/2)! 2^{j/2}} \Biggl(\frac{1}{2} \sum_{k < \frac{q\log(q/H_n)}{H_n}} |a_k|^2 \Biggr)^{j/2} \Biggr)^2 , \nonumber
\end{eqnarray}
where $f(k)$ is an extended Rademacher random multiplicative function. As we noted in section \ref{momentsmanipsec}, the $a_k$ are bounded in absolute value by 1, so the error term on the third line is $\ll \frac{Q^{0.26}}{Q^{0.99}} \left( (\frac{q\log(q/H_n)}{H_n})^j + j^{j/2} (\frac{q\log(q/H_n)}{H_n})^{j/2} \right)^2 \ll \frac{Q^{0.78}}{Q^{0.99}}$, which is negligible. 

Finally we apply the second part of Probability Result \ref{rmfimport} to handle the expectation in the above display. Since the $a_k$ are real valued, bounded in absolute value by 1, and we only need to establish \eqref{neededafterreg} for $1 \leq j \leq \min\{\frac{\log^{1/4}(Q/H_{n}(Q))}{50}, \frac{\log Q}{4\log(Q/H_{n}(Q))}\} \leq \frac{\log^{1/4}(Q/H_{n}(Q))}{50}$, we see all the conditions of Probability Result \ref{rmfimport} are satisfied. Recall once more that $q/H_{n}(q) = Q/H_{n}(Q)$ is constant for all $Q \leq q \leq 2Q$. So our expectation is
\begin{eqnarray}
& \ll & \left(\frac{Q}{H_{n}(Q)} \right)^{j - 1/15j} \log^{j}(Q/H_{n}(Q)) = \left(\frac{Q}{4 H_{n}(Q)} \right)^{j} e^{- \frac{\log(Q/H_{n}(Q))}{15j}} (4\log(Q/H_{n}(Q)))^j \nonumber \\
& \leq & \left(\frac{Q}{4 H_{n}(Q)} \right)^{j} e^{- \frac{10\log^{3/4}(Q/H_{n}(Q))}{3}} (4\log(Q/H_{n}(Q)))^{\frac{\log^{1/4}(Q/H_{n}(Q))}{50}} , \nonumber
\end{eqnarray}
which is more than good enough to imply \eqref{neededafterreg}.
\qed

% SECTION 7 %%%%%%%%%%%%%%%%%%%%%%%%%%%%%%
\section{Proof of Theorem 4}
The proof of Theorem 4 is very similar to, but simpler than, the proof of Theorem 3.

Recall the reductions from section \ref{fourierred}, and that our target distribution in Theorem 4 is $Z_1 + iZ_2$ with $Z_1, Z_2$ independent $N(0,1/2)$ random variables, having moments $\E (Z_1 + iZ_2)^j (\overline{Z_1 + iZ_2})^k = \frac{1}{\pi} \int_{-\infty}^{\infty} \int_{-\infty}^{\infty} (z_1 + iz_2)^j (z_1 - iz_2)^k e^{-z_1^2 - z_2^2} dz_1 dz_2 = k! \textbf{1}_{j=k}$. (This is easy to check after rewriting the double integral in polar coordinates.) Then it will suffice to show the existence of sets $\mathcal{G}_{q,H}$ of characters mod $q$, with the sizes claimed in the theorem, such that for any choice of $\chi \in \mathcal{G}_{q,H}$ and all $j, k \geq 0$ we have
\begin{eqnarray}
&& \int_{0}^{1} \Biggl( \frac{2\sqrt{q}}{\pi \sqrt{H}} \sum_{1 \leq m < (q/H)\log(q/H)} \frac{\overline{\chi}(m) \sin(\pi m H/q)}{m} \cos(2\pi m\theta) \Biggr)^{j} \cdot \nonumber \\
&& \cdot \Biggl( \overline{\frac{2\sqrt{q}}{\pi \sqrt{H}} \sum_{1 \leq m < (q/H)\log(q/H)} \frac{\overline{\chi}(m) \sin(\pi m H/q)}{m} \cos(2\pi m\theta)} \Biggr)^{k} d\theta \rightarrow k! \textbf{1}_{j=k} \;\;\;\;\; \text{as} \; q \rightarrow \infty . \nonumber
\end{eqnarray}
(As for Theorem 3, we actually need to show this with $\cos(2\pi m\theta)$ replaced by $\sin(2\pi m\theta)$ as well, but that case will be exactly similar so we shall not discuss it further.)

If we set $a_m = a_{m,q,H} := \frac{q \sin(\pi m H/q)}{\pi H m}$, as in section \ref{momentsmanipsec}, then we can rewrite our goal as being that for any choice of $\chi \in \mathcal{G}_{q,H}$ and all $j, k \geq 0$, we have
\begin{eqnarray}
&& \left(\frac{4H}{q} \right)^{(j+k)/2} \Biggl| \int_{0}^{1} \Biggl( \sum_{m < (q/H)\log(q/H)} a_m \overline{\chi}(m) \cos(2\pi m\theta) \Biggr)^{j} \Biggl( \overline{\sum_{m < (q/H)\log(q/H)} a_m \overline{\chi}(m) \cos(2\pi m\theta)} \Biggr)^{k} - \nonumber \\
&& \;\;\;\;\;\;\;\;\;\;\;\;\;\;\;\;\;\;\;\; - k! \textbf{1}_{j=k} \Biggl(\frac{1}{2} \sum_{m < \frac{q\log(q/H)}{H}} |a_m|^2 \Biggr)^{(j+k)/2} \Biggr| \rightarrow 0 \;\;\; \text{as} \; q \rightarrow \infty . \nonumber
\end{eqnarray}
To establish this, it will suffice to show that for all $0 \leq j, k \leq \min\{\frac{\log^{1/4}(q/H(q))}{50}, \frac{\log q}{4\log(q/H(q))}\}$, say, we have
\begin{eqnarray}\label{neededcomplex}
&& \frac{1}{q-1} \sum_{\chi \; \text{mod} \; q} \Biggl| \int_{0}^{1} \Biggl( \sum_{m < \frac{q\log(q/H)}{H}} a_m \overline{\chi}(m) \cos(2\pi m\theta) \Biggr)^{j} \Biggl( \overline{\sum_{m < \frac{q\log(q/H)}{H}} a_m \overline{\chi}(m) \cos(2\pi m\theta)} \Biggr)^{k} - \nonumber \\
&& \;\;\;\;\;\;\;\;\;\;\;\;\;\;\;\;\;\;\;\; - k! \textbf{1}_{j=k} \Biggl(\frac{1}{2} \sum_{m < \frac{q\log(q/H)}{H}} |a_m|^2 \Biggr)^{(j+k)/2} \Biggr|^2 \ll \left(\frac{q}{4H} \right)^{j+k} e^{-2\log^{3/4}(q/H)} .
\end{eqnarray}
For if we have \eqref{neededcomplex}, then the proportion of $\chi$ mod $q$ for which
\begin{eqnarray}
&& \left(\frac{4H}{q} \right)^{j+k} \Biggl| \int_{0}^{1} \Biggl( \sum_{m < (q/H)\log(q/H)} a_m \overline{\chi}(m) \cos(2\pi m\theta) \Biggr)^{j} \Biggl( \overline{\sum_{m < (q/H)\log(q/H)} a_m \overline{\chi}(m) \cos(2\pi m\theta)} \Biggr)^{k} - \nonumber \\
&& \;\;\;\;\;\;\;\;\;\;\;\;\;\;\;\;\;\;\;\; - k! \textbf{1}_{j=k} \Biggl(\frac{1}{2} \sum_{m < \frac{q\log(q/H)}{H}} |a_m|^2 \Biggr)^{(j+k)/2} \Biggr|^2 \gg e^{-0.1\log^{3/4}(q/H)} \nonumber
\end{eqnarray}
must be $\ll e^{-1.9\log^{3/4}(q/H)}$, and so the proportion for which this holds for some pair of $0 \leq j, k \leq \min\{\frac{\log^{1/4}(q/H)}{50}, \frac{\log q}{4\log(q/H)}\}$ must be $\ll \log^{1/2}(q/H) e^{-1.9\log^{3/4}(q/H)} \ll e^{-1.8\log^{3/4}(q/H)}$. Excluding any such characters mod $q$, our remaining set $\mathcal{G}_{q,H}$ of ``good'' characters will satisfy $\#\mathcal{G}_{q,H} \geq q(1 - O(e^{-1.8\log^{3/4}(q/H)}))$, which is more than good enough for Theorem 4. And under the hypotheses of the theorem, we have $e^{-0.1\log^{3/4}(q/H)} \rightarrow 0$ as well as $\min\{\frac{\log^{1/4}(q/H(q))}{50}, \frac{\log q}{4\log(q/H(q))}\} \rightarrow \infty$ as $q \rightarrow \infty$, so for any fixed $j,k$ and for $\chi \in \mathcal{G}_{q,H}$ the moment will tend to the desired Gaussian moment.

Now it only remains to verify \eqref{neededcomplex}. But expanding the square on the left hand side there, using multiplicativity of $\chi$ and the condition that $j, k \leq \frac{\log q}{4\log(q/H)}$, we see the resulting expression only involves $\chi$ and $\overline{\chi}$ applied to numbers that are $\leq ((q/H)\log(q/H))^{j+k} \leq q^{0.51}$, say (for large enough $q$). Since this is $< q$, the orthogonality of Dirichlet characters implies that the left hand side in \eqref{neededcomplex} is exactly equal to
\begin{eqnarray}
&& \E \Biggl| \int_{0}^{1} \Biggl( \sum_{m < \frac{q\log(q/H)}{H}} a_m f(m) \cos(2\pi m\theta) \Biggr)^{j} \Biggl( \overline{\sum_{m < \frac{q\log(q/H)}{H}} a_m f(m) \cos(2\pi m\theta)} \Biggr)^{k} d\theta - \nonumber \\
&& \;\;\;\;\;\;\;\;\;\;\;\;\;\;\;\;\;\;\;\; - k! \textbf{1}_{j=k} \Biggl(\frac{1}{2} \sum_{m < \frac{q\log(q/H)}{H}} |a_m|^2 \Biggr)^{(j+k)/2} \Biggr|^2 , \nonumber
\end{eqnarray}
where $f(m)$ is a Steinhaus random multiplicative function. The desired bound now follows immediately from Probability Result \ref{rmfsteinhausimport} and a small computation.
\qed

\vspace{12pt}
\noindent {\em Acknowledgements.}
The author would like to thank K. Soundararajan and Max Xu for sharing a draft of their forthcoming paper~\cite{soundxu}, and Max Xu for some helpful comments.


\begin{thebibliography}{99}

\bibitem{bennishrodg} J. Benatar, A. Nishry, B. Rodgers. Moments of polynomials with random multiplicative coefficients. Preprint available online at \url{https://arxiv.org/abs/2012.15507}

\bibitem{bobgoldgrankou} J. Bober, L. Goldmakher, A. Granville, D. Koukoulopoulos. The frequency and the structure of large character sums. {\em J. Eur. Math. Soc. (JEMS)}, \textbf{20}, no. 7, pp 1759-1818. 2018

\bibitem{burkholderaoms} D. L. Burkholder. Martingale transforms. {\em Ann. Math. Statist.}, \textbf{37}, pp 1494-1504. 1966

\bibitem{chatsound} S. Chatterjee, K. Soundararajan. Random multiplicative functions in short intervals. {\em Int. Math. Res. Not.}, pp 479-492. 2012

\bibitem{davenporterdos} H. Davenport, P. Erd\H{o}s. The distribution of quadratic and higher residues. {\em Publ. Math. Debrecen}, \textbf{2}, 252-265. 1952 



\bibitem{gransoundlcs} A. Granville, K. Soundararajan. Large character sums. {\em J. Amer. Math. Soc.}, \textbf{14}, no. 2, pp 365-397. 2001

\bibitem{gransoundpv} A. Granville, K. Soundararajan. Large character sums: pretentious characters and the P\'olya--Vinogradov theorem. {\em J. Amer. Math. Soc.}, \textbf{20}, no. 2, pp 357-384. 2007

\bibitem{gut} A. Gut. {\em Probability: A Graduate Course.} Second edition, published by Springer Texts in Statistics. 2013



\bibitem{harperlimits} A. J. Harper. On the limit distributions of some sums of a random multiplicative function. {\em Journal f\"{u}r die reine und angewandte Mathematik}, \textbf{678}, pp 95-124. 2013

\bibitem{harperrmflow} A. J. Harper. Moments of random multiplicative functions, I: Low moments, better than squareroot cancellation, and critical multiplicative chaos. {\em Forum of Mathematics, Pi}, \textbf{8}, e1, 95pp. 2020

\bibitem{hughesrudnick} C. P. Hughes, Z. Rudnick. On the distribution of lattice points in thin annuli. {\em Int. Math. Res. Not.}, \textbf{2004}, no. 13, pp 637-658. 2004

\bibitem{hussain} A. Hussain. The limiting distribution of character sums. Preprint available online at \url{https://arxiv.org/abs/2010.06967}

\bibitem{kalmynin} A. B. Kalmynin. Large values of short character sums. {\em J. Number Theory}, \textbf{198}, pp 200-210. 2019

\bibitem{lamzouridcs} Y. Lamzouri. The distribution of short character sums. {\em Math. Proc. Cambridge Philos. Soc.}, \textbf{155}, no. 2, pp 207-218. 2013

\bibitem{makzah} K.-H. Mak, A. Zaharescu. The distribution of values of short hybrid exponential sums on curves over finite fields. {\em Math. Res. Lett.}, \textbf{18}, no. 1, pp 155-174. 2011


\bibitem{mvcharmean} H. L. Montgomery, R. C. Vaughan. Mean values of character sums. {\em Canadian J. Math.}, \textbf{31}, no. 3, pp 476-487. 1979 

\bibitem{mv} H. L. Montgomery, R. C. Vaughan. {\em Multiplicative Number Theory I: Classical Theory.} First edition, published by Cambridge University Press. 2007

\bibitem{najnudelconsec} J. Najnudel. On consecutive values of random completely multiplicative functions. {\em Electron. J. Probab.}, \textbf{25}, Paper No. 59, 28 pp. 2020

\bibitem{perret} C. Perret-Gentil. Gaussian distribution of short sums of trace functions over finite fields. {\em Math. Proc. Cambridge Philos. Soc.}, \textbf{163}, no. 3, pp 385-422. 2017 

\bibitem{salzyg} R. Salem, A. Zygmund. Some properties of trigonometric series whose terms have random signs. {\em Acta Math.}, \textbf{91}, no. 1, pp 245-301. 1954

\bibitem{soundxu} K. Soundararajan, M. W. Xu. Central limit theorems for random multiplicative functions. In preparation.

\bibitem{vaughanwooley} R. C. Vaughan, T. D. Wooley. On a certain nonary cubic form and related equations. {\em Duke Math. J.}, \textbf{80}, no. 3, pp 669-735. 1995


\end{thebibliography}
\end{document}